\documentclass[10pt]{article}
\usepackage{amssymb,amsmath,amsthm}
\usepackage{indentfirst}
\usepackage{exscale}
\usepackage{relsize}
\usepackage{subfigure}
\usepackage{graphicx}
\newcommand{\Z}{\mathbb{Z}}
\newcommand{\R}{\mathbb{R}}
\theoremstyle{definition}

\theoremstyle{remark}

\numberwithin{equation}{section}
\allowdisplaybreaks
 \UseRawInputEncoding
\begin{document}
\title{\Large\bf{ Existence and multiplicity of nontrivial solutions for a $(p,q)$-Laplacian system on locally finite graphs}}
\date{}
\author {Ping Yang$^1$, \ Xingyong Zhang$^{1,2}$\footnote{Corresponding author, E-mail address: zhangxingyong1@163.com}\\
{\footnotesize $^1$Faculty of Science, Kunming University of Science and Technology,}\\
 {\footnotesize Kunming, Yunnan, 650500, P.R. China.}\\
{\footnotesize $^{2}$Research Center for Mathematics and Interdisciplinary Sciences, Kunming University of Science and Technology,}\\
 {\footnotesize Kunming, Yunnan, 650500, P.R. China.}\\
      }

 \date{}
 \maketitle

\begin{center}
\begin{minipage}{15cm}
\par
\small  {\bf Abstract:} We generalize two embedding theorems and investigate the existence and multiplicity of nontrivial solutions for a $(p,q)$-Laplacian coupled system with perturbations
and two parameters $\lambda_1$ and $\lambda_2$ on locally finite graph. By using the Ekeland's variational principle, we obtain that system has at least
one nontrivial solution when the nonlinear term satisfies
 the sub-$(p,q)$ conditions. We also obtain a necessary condition for the existence of  semi-trivial solutions to the system. Moreover,
by using the mountain pass theorem and Ekeland's variational principle,
 we obtain that  system has at least one  solution of positive energy and one  solution of negative energy when the nonlinear term satisfies
 the super-$(p,q)$ conditions which is weaker than the well-known Ambrosetti-Rabinowitz condition. Especially, in all of the results, we present the concrete ranges of the parameters $\lambda_1$ and $\lambda_2$.
\end{minipage}
 \end{center}

\section{Introduction }
\setcounter{equation}{0}
Let $G=(V, E)$ be a locally finite and connected graph, where $V$ denotes the vertex set and $E$ denotes the edge set. We say that $(V,E)$ is a locally finite graph if for any $x \in V$ there are only finite edges $xy \in E$.
 Moreover, We say that $(V,E)$ is a connected graph if any two vertices $x$ and $y$ can be connected via finite edges.
 For any edge $xy\in E$, assume that its weight $\omega_{xy}>0$ and $\omega_{xy}=\omega_{yx}$.
 For any $x\in V$, its degree is defined as $\mbox{deg}(x)=\sum_{y\thicksim x}\omega_{xy}$, where we denote $y\thicksim x$ if there exists $y \in V$ such that edge $xy\in E$.
 The distance of two vertices $x,y$, denoted by $dist(x,y)$, is defined as the minimal number of edges which connect $x,y$.
 Let $\mu:V\to \R^+$ be a finite measure, $\mu(x)\geq \mu_0>0$,
 and $C(V)$ be the set of all real functions on $V$. Define $\Delta: C(V)\to C(V)$ as
\begin{eqnarray}
\label{eq8}
\Delta u(x)=\frac{1}{\mu(x)}\sum_{y\thicksim x}\omega_{xy}(u(y)-u(x)).
\end{eqnarray}
The associated gradient form is
\begin{eqnarray}
\label{eq9}
\Gamma(u,v)(x)=\frac{1}{2\mu(x)}\sum\limits_{y\thicksim x}w_{xy}(u(y)-u(x))(v(y)-v(x)).
\end{eqnarray}
Write $\Gamma(u)=\Gamma(u,u)$. We denote the length of the gradient is
\begin{eqnarray}
\label{eq10}
|\nabla u|(x)=\sqrt{\Gamma(u)(x)}=\left(\frac{1}{2\mu(x)}\sum\limits_{y\thicksim x}w_{xy}(u(y)-u(x))^2\right)^{\frac{1}{2}}.
\end{eqnarray}
We can obtain that the gradient has the following properties:
\begin{eqnarray}
(1)&&\Gamma(u_1+u_2,v)(x)   =  ~~\Gamma(u_1,v)(x)+\Gamma(u_2,v)(x);\nonumber\\
(2)&&\Gamma(u,v_1+v_2)(x)   =  ~~\Gamma(u,v_1)(x)+\Gamma(u,v_2)(x);\nonumber\\
(3)&&\Gamma(\theta u,v)(x)  =  ~~\Gamma(u,\theta v)(x)=\theta\Gamma(u,v)(x)\mbox{  for all }\theta \in \R;\nonumber\\
(4)&&\Gamma(u,v)        \le ~~|\nabla u|\cdot|\nabla v|;\label{G1}\\
(5)&&\big||\nabla u_k|-|\nabla u|\big|\le ~~|\nabla(u_k-u)|.\label{G2}
\end{eqnarray}
In fact, (1)-(3) obviously hold by the definition of $\Gamma(u,v)$.  By (\ref{eq9}), we have
\begin{eqnarray*}
\Gamma(u,v)& = & \frac{1}{2\mu(x)}\sum\limits_{y\thicksim x}w_{xy}(u(y)-u(x))(v(y)-v(x))\\
              &\le& \left(\frac{1}{2\mu(x)}\sum\limits_{y\thicksim x}w_{xy}(u(y)-u(x))^2\right)^{\frac{1}{2}}\left(\frac{1}{2\mu(x)}\sum\limits_{y\thicksim x}w_{xy}(v(y)-v(x))^2\right)^{\frac{1}{2}}\\
              & = & |\nabla u|\cdot|\nabla v|
\end{eqnarray*}
Moreover, by (\ref{eq10}) and (\ref{G1}), we have
\begin{eqnarray*}
        |\nabla(u_k-u)|^2
&  =  & \Gamma(u_k-u,u_k-u)\\
&  =  & \frac{1}{2\mu(x)}\sum\limits_{y\thicksim x}w_{xy}\left((u_k(y)-u(y))-(u_k(x)-u(x))\right)^2\\
&  =  & \frac{1}{2\mu(x)}\sum\limits_{y\thicksim x}w_{xy}\left((u_k(y)-u_k(x))-(u(y)-u(x))\right)^2\\
&  =  & \frac{1}{2\mu(x)}\sum\limits_{y\thicksim x}w_{xy}[(u_k(y)-u_k(x))^2+(u(y)-u(x))^2-2(u(y)-u(x))(u_k(y)-u_k(x))]\\
&  =  & |\nabla u_k|^2+|\nabla u|^2-2\Gamma(u,u_k)\\
& \ge & |\nabla u_k|^2+|\nabla u|^2-2|\nabla u_k||\nabla u|\\
&  =  & (|\nabla u_k|-|\nabla u|)^2.
\end{eqnarray*}
Hence, (\ref{G2}) holds.
\par
For any $p> 1$, we define $ \Delta_{p}:C(V)\rightarrow C(V) $ as follows:
\begin{equation}\label{c1}
  \Delta_{p}(u)(x) = \frac{1}{2\mu(x)} \sum_{y\sim x}  \left(|\nabla u|^{p-2} (y)+ |\nabla u|^{p-2} (x) \right) \omega_{xy} (u(y)-u(x)).
  \end{equation}
Let $C_c(V):=\{u:V \to \R \big| suppu\subset V \}$. Then for any function $\phi \in C_c(V)$,
\begin{equation}\label{p1}
\int_V \Delta_{p}u\phi d\mu=-\int_V |\nabla u|^{p-2}\Gamma(u,\phi)d\mu.
 \end{equation}
For any function $u:V\to\R$, we denote
\begin{eqnarray}
\label{eq12}
\int\limits_V u(x) d\mu=\sum\limits_{x\in V}u(x)\mu(x).
\end{eqnarray}
Define $L^r(V)=\{u:V \to \R \big|\int_V |u|^r d\mu <+\infty\}$ ($1\leq r<+\infty$) with the norm defined by
$$
\|u\|_{L^r(V)}=
\left(\int_V|u(x)|^r d \mu\right)^{\frac{1}{r}}.
$$
Then  $(L^r(V),\|\cdot\|_{L^r(V)})$ is a Banach space  ($1\le r<+\infty$).
Define $L^\infty(V)=\left\{u:V \to \R \Big| \sup_{x\in V}|u(x)|<+\infty\right\}$  with the norm defined by
$$
\|u\|_{\infty}=\sup_{x\in V}|u(x)|.
$$
For more details, one can see \cite{Grigor'yan 2017,Yamabe 2016}.
\par
Consider the following $p$-Laplacian equation on a locally finite graph  $G=(V,E)$,
\begin{eqnarray}
\label{aa1}
-\Delta_p u+h(x)|u|^{p-2}u=f(x,u),\ \ x\in V,
\end{eqnarray}
where $p>1$, $h:V\to \R$ and $f:V\times \R\to \R$.
\par
In recent years, the existence and multiplicity of nontrivial solutions to (\ref{aa1}) have attracted some attentions
(for example, see \cite{2 order positive solution,zhang and lin 2019, zhang and chang 2021,shaomengqiu 2023, man 2020, liuyang 2022, Grigor'yan 2017, Yamabe 2016}).
 In \cite{2 order positive solution}, Zhang investigated (\ref{aa1}) with $p=2$ and $f(x,u)=|u|^{s-2}u$
for all $x\in V$, where $s>2$. He obtained that equation (\ref{aa1}) has a positive solution by using the mountain pass theorem.
 In \cite{ zhang and lin 2019}, Zhang and Lin studied (\ref{aa1}) with $f(x,u)=g(x)|u|^{r-2}u$ for all $x\in V$, where $g:V\to \R$ and $r>p>2$. They obtained that equation (\ref{aa1}) has a positive solution.
 In \cite{zhang and chang 2021}, by using the variational principles and Fatou's lemma, Chang and Zhang obtained the equation (\ref{aa1}) has a solution when $f(x,u)$ is Lipschitz continuous in $u$.
In \cite{shaomengqiu 2023}, Shao investigated (\ref{aa1}) with $f(x,u)=g(x,u)+e(x)$. When $\|e\|_{L^{\frac{p}{p-1}}(V)}$ small enough, $g(x,u)$ satisfies sub-$(p-1)$-linear growth condition at origin and $|g(x,u)|<C\left(1+{|u|}^{q-1}\right)$ for all $x\in V, \; \mbox{where} \; q>p\geq 2$, Shao obtained the equation (\ref{aa1}) has one nontrivial solution of positive energy and another nontrivial solution of negative energy by using the mountain pass theorem and Ekeland's variational principle.
In \cite{man 2020}, Man investigated (\ref{aa1}) with $p=2$ and $h$ replaced by a constant $\alpha$. When $\alpha$ small enough and nonlinear term $f(x,u)$ satisfies super-$(r-1)$-linear growth condition at origin, where $r>2$ and some additional assumptions, he obtained that equation (\ref{aa1}) has a positive solution by using the mountain pass theorem.
In \cite{liuyang 2022}, Liu investigated (\ref{aa1}) with $p=2$  and Dirichlet boundary condition, where $f(x,u)=|u|^{r-2}u+\epsilon e(x)$, where $r>2$, $\epsilon>0$ and $e(x)>0$. When $\epsilon$ small enough, he obtained that the equation has two positive solutions by using the mountain pass theorem and Ekeland's variational principle.
Especially, in \cite{Grigor'yan 2017}, Grigor'yan-Lin-Yang considered (\ref{aa1}) with $p=2$.
They assumed that the measure $\mu(x)\geq \mu_{\min}>0$ for all $x \in V$, where $\mu_{\min}=\min_{x\in V} \mu(x)$, and $h$ and $f$ satisfy the following conditions:\\
$(K_1)$ \; there exists a constant $h_0>0$ such that $h(x)\geq h_0$ for all $x\in V$;\\
$(K_2)$ \; $\frac{1}{h} \in L^1(V)$;\\
$(S_1)$ \; $f(x,s)$ is continuous in $s$, $f(x,0)=0$, and for any fixed $M>0$, there exists a constant $A_M$ such that $\max_{s \in [0,M]}f(x,s)\le A_M $ for all $x \in V$;\\
$(S_2)$ \; $\limsup_{s\rightarrow 0^+} \frac{2F(x,s)}{s^2}<\lambda_1=\inf_{\int_V u^2 d\mu=1} \int_V (|\nabla u|^2+hu^2 )d\mu;$\\
$(S_3)$ \; there exists a constant $\theta>2$ such that for all $x\in V$ and $s>0$,
$$
0<\theta F(x,s)=\theta\int_0^s f(x,t) dt\le sf(x,s).
$$
(The $(S_3)$ condition is usually called as Ambrisetti-Rabinowitz condition ((AR)-condition for short).)\\
Then equation (\ref{aa1}) with $p=2$ has a strictly positive solution.
 Moreover, they also investigated the following equation with perturbation:
\begin{eqnarray}
\label{aa2}
 -\Delta u+hu=f(x,u)+ \epsilon e(x),\ \ x\in V,
\end{eqnarray}
 where $e\geq0$ for all $x\in V$ ($e\not\equiv 0$). They obtained that there exists a constant $\epsilon_0>0$ such that for any $0<\epsilon<\epsilon_0$,  (\ref{aa2}) has at least two distinct strictly positive solutions under the above assumptions. When $(K_2)$ is replaced by the following condition: \\
 $(K_2')$ \; $h(x)\rightarrow +\infty$ as dist$(x,x_0)\rightarrow +\infty$ for some fixed $x_0 \in V$,\\
  and $(S_1)$ is replaced by the following condition:\\
  $(S_1')$ \; $f(x,0)=0,\;f(x,s)>0$ for all $x\in V$ and all $s>0$, and there exists a constant $L>0$ such that
$$
|f(x,s)-f(x,t)|\le L|s-t|\ \ \mbox{for all}\; x\in V\mbox{ and all }(s,t) \in \R^2.
$$
They obtained that (\ref{aa2}) has a strictly positive solution.
\par
In this paper,  inspired by \cite{Grigor'yan 2017,Yamabe 2016}
we consider the following $(p,q)$-Laplacian coupled system with perturbation terms and two parameters on  a locally finite graph  $G=(V,E)$:
\begin{eqnarray}
\label{eq1}
 \begin{cases}
   -\Delta_p u+h_1(x)|u|^{p-2}u=F_u(x,u,v)+\lambda_1 e_1(x),\;\;\;\;\hfill x\in V,\\
   -\Delta_q v+h_2(x)|v|^{q-2}v=F_v(x,u,v)+\lambda_2 e_2(x),\;\;\;\;\hfill x\in V,\\
 \end{cases}
\end{eqnarray}
where $-\Delta_p$ and $-\Delta_q$ are defined by (\ref{c1}) with $p\ge2$ and $q\ge2$, $F:V\times \R^2 \to \R$,
$e_1\in L^{\frac{p}{p-1}}(V),\;e_2\in L^{\frac{q}{q-1}}(V),\;e_1(x),e_2(x)\not\equiv 0$ and $\lambda_1,\lambda_2> 0$.
\par
If $(u,v)$ is a solution of system (\ref{eq1}) and $(u,v)\not=(0,0)$, then we call that $(u,v)$ is a nontrivial solution of system (\ref{eq1}).
Furthermore, if $(u,v)$ is a nontrivial solution of system (\ref{eq1}),  $(u,v)=(u,0)$ or $(u,v)=(0,v)$, then
 we call that $(u,v)$ is a semi-trivial solution of system (\ref{eq1}).
We obtain the following results.
\vskip2mm
\noindent
$\bullet$ {\bf (I) The sub-$(p,q)$-linear case}
\vskip2mm
\noindent
{\bf Theorem 1.1.} {\it Assume that the following conditions hold:\\
$(H_1)$ \; there exists a constant $h_0>0$ such that $ h_i(x)\ge h_0>0$ for all $x\in V$, $i=1,2$;\\
$(H_2)$ \; $h_i(x)\rightarrow \infty$ as dist$(x,x_0)\rightarrow \infty$ for some fixed $x_0,\; i=1,2$;\\
$(F_0)$\; $F(x,s,t)$ is continuously differentiable in $(s,t)\in \R^2$ for all $x\in V$,  and there exists a  function $a\in C(\R^+,\R^+)$ and a function $b:V\to \R^+$ with $b\in L^1(V)$ such that
$$
|F_s(x,s,t)|\le a(|(s,t)|) b(x), |F_t(x,s,t)|\le a(|(s,t)|) b(x), |F(x,s,t)|\le a(|(s,t)|) b(x),
$$
for all $x\in V$ and all $(s,t)\in \R^2$;\\
$(F_1)$ \;  $F(x,0,0)=0$, and there exists $f_i,g_i:V\to \R^+,\;i=1,2$, $g_1\in L^{\frac{p}{p-1}}(V)$
 and $g_2\in L^{\frac{q}{q-1}}(V)$ with $ \|f_1\|_{\infty}<\min\{\frac{h_0}{2},\frac{ph_0}{q(p-1)}\}$
  and $\|f_2\|_{\infty}<h_0-\frac{q(p-1)}{p}\|f_1\|_{\infty}$ such that
$$
|F_s(x,s,t)|\le f_1(x)(|s|^{p -1}+|t|^{\frac{pq-q}{p}})+g_1(x),\;
|F_t(x,s,t)|\le f_2(x)(|s|^p+|t|^{q-1})+g_2(x),
$$
for all $x\in V$ and all $(s,t)\in \R^2$;\\
$(F_2)$ one of the following conditions holds:
\par
(i) there exists $\beta_1>1$ and $K_1:V\to \R$ such that $K_1(x_1)>0$ for some $x_1 \in V$ with $e_1(x_1)>0$ and
$$
 F(x,s,0)\ge -K_1(x)|s|^{\beta_1} \ \ \mbox{for all } s\in \R \mbox{ and all  } x\in V;
$$
\par
(ii) there exists $\beta_2>1$ and $K_2:V\to \R$ such that $K_2(x_2)>0$ for some $x_2 \in V$ with $e_2(x_2)>0$ and
$$
 F(x,0,t)\ge -K_2(x)|t|^{\beta_2} \ \ \mbox{for all } t\in \R \mbox{ and all  } x\in V.
$$
 Then for each pair $(\lambda_1,\lambda_2)\in (0,+\infty)\times(0,+\infty)$, system (\ref{eq1})
 has at least  one nontrivial solution $(u_{\lambda \star},v_{\lambda \star})$.
 Furthermore, we show the necessary conditions for the existence of the non-semi-trivial solutions to the system (\ref{eq1}).
 If $(u_{\lambda \star},v_{\lambda \star})=(u_{\lambda \star},0)$, then
$$
\|u_{\lambda \star}\|_{\infty}\le {\mu_0}^{-\frac{1}{p}}\left(\frac{\lambda_1\|e_1\|_{L^{\frac{p}{p-1}}(V)}+\|g_1\|_{L^{\frac{p}{p-1}}(V)}}
{h_0-\|f_1\|_{\infty}}\right)^{\frac{1}{p-1}}.
$$
If $(u_{\lambda \star},v_{\lambda \star})=(0,v_{\lambda \star})$, then}
$$
\|v_{\lambda \star}\|_{\infty}\le {\mu_0}^{-\frac{1}{q}}\left(\frac{\lambda_2\|e_2\|_{L^{\frac{q}{q-1}}(V)}+\|g_2\|_{L^{\frac{q}{q-1}}(V)}}
{h_0-\|f_2\|_{\infty}}\right)^{\frac{1}{q-1}}.
$$

\vskip2mm
\noindent
{\bf Theorem 1.2.} {\it  Assume that $(H_1)$, $(H_2)$, $(F_0)$, $(F_2)$ and  the following condition hold:\\
$(F_1')$ \;  $F(x,0,0)=0$, and there exists $f_i,g_i:V\to \R^+,\;i=1,2$,  $g_1\in L^{\frac{p}{p-1}}(V)$ and $g_2\in L^{\frac{q}{q-1}}(V)$
with $ \|f_1\|_{\infty}<\min\{\frac{h_0}{2},\frac{qh_0}{p(q-1)}\}$ and $\|f_2\|_{\infty}<h_0-\frac{p(q-1)}{q}\|f_2\|_{\infty}$
such that
$$
|F_s(x,s,t)|\le f_2(x)(|t|^q+|s|^{p-1})+g_1(x), \;
|F_t(x,s,t)|\le f_1(x)(|t|^{q -1}+|s|^{\frac{qp-p}{q}})+g_2(x),
$$
for all $x\in V$ and all $(s,t)\in \R^2$.\\
 Then for each pair $(\lambda_1,\lambda_2)\in (0,+\infty)\times(0,+\infty)$, system (\ref{eq1}) has at least
 one nontrivial solution $(u_{\lambda \star},v_{\lambda \star})$.
 Furthermore, we show the necessary conditions for the existence of the non-semi-trivial solutions to the system (\ref{eq1}).
 If $(u_{\lambda \star},v_{\lambda \star})=(u_{\lambda \star},0)$, then
$$
 \|u_{\lambda \star}\|_{\infty}\le {\mu_0}^{-\frac{1}{p}}\left(\frac{\lambda_1\|e_1\|_{L^{\frac{p}{p-1}}(V)}+\|g_1\|_{L^{\frac{p}{p-1}}(V)}}
{h_0-\|f_2\|_{\infty}}\right)^{\frac{1}{p-1}}.
$$
If $(u_{\lambda \star},v_{\lambda \star})=(0,v_{\lambda \star})$, then}
$$
\|v_{\lambda \star}\|_{\infty}\le {\mu_0}^{-\frac{1}{q}}\left(\frac{\lambda_2\|e_2\|_{L^{\frac{q}{q-1}}(V)}+\|g_2\|_{L^{\frac{q}{q-1}}(V)}}
{h_0-\|f_1\|_{\infty}}\right)^{\frac{1}{q-1}}.
$$

\vskip2mm
\noindent
$\bullet$  {\bf (II) The super-$(p,q)$-linear case}
\vskip2mm
 \noindent
 {\bf Theorem 1.3.} {\it Let $\lambda_1=\lambda_2=\lambda$. Assume $(H_1)$,  $(F_0)$ and the following conditions hold: \\
$(H_2')$ \; for any given constant $B>0$, $\sum_{x\in A_i} \mu(x)<\infty$, where $A_i=\{x\in V|h_i(x)\le B\}$, $i=1,2$;\\
$(C_1)$ \; $F(x,0,0)=0$ for all $x\in V$, and there exists a constant $l_0>0$ such that
$$
|F_s(x,s,t)|\le \frac{h_0}{q+1}(|s|^{p -1}+|t|^{\frac{pq-q}{p}}),\;|F_t(x,s,t)|\le \frac{h_0}{q+1}(|s|^p+|t|^{q-1}),
$$
for and all $x\in V$  and all $(s,t)\in \R^2$ with $|(s,t)|< l_0$; \\
$(C_2)$ \;  there exists $l_1>0$ such that $F(x_3,s,s)\geq M(s^p+s^q)\;\mbox{for some } x_3\in \{x\in V \big|e_1(x)+e_2(x)>0\}$ and all $s\in \R$ with
$s> l_1$, where
$$
M> \max\left\{\frac{D_1+\mu(x_3)h_1(x_3)}{p\mu(x_3)},
\frac{D_2+\mu(x_3)h_2(x_3)}{q\mu(x_3)}\right\},
$$
 $D_1=\left(\frac{deg(x_3)}{2}\right)^{\frac{p}{2}}\left(\sum_{x\thicksim x_3}\left(\frac{1}{\mu(x)}\right)^{\frac{p}{2}-1}+\frac{1}{\mu(x_3)^{\frac{p}{2}-1}}\right)$ and  $D_2=\left(\frac{deg(x_3)}{2}\right)^{\frac{q}{2}}\left(\sum_{x\thicksim x_3}\left(\frac{1}{\mu(x)}\right)^{\frac{q}{2}-1}+\frac{1}{\mu(x_3)^{\frac{q}{2}-1}}\right)$;\\
$(C_3)$ \; there exists a constant $\nu>\max\{p,q\}$ and $0\le A<\min\left\{\dfrac{\nu}{p}-1,\dfrac{\nu}{q}-1\right\}h_0$ such that
$$
\nu F(x,s,t)-F_s(x,s,t)s-F_t(x,s,t)t\le A( |s|^p+|t|^q)  \; \; \mbox{for all}\;x \in V.
$$
Then for each $\lambda$ satisfying
\begin{eqnarray}
\label{pp1}
0<\lambda<\lambda_0=\frac{\min\{1,q-1\}}{2^{\max\{p,q\}-1}(pq+p)\max\left\{{h_0}^{-\frac{1}{p}}\|e_1\|_{L^{\frac{p}{p-1}}(V)}
,{h_0}^{-\frac{1}{q}}\|e_2\|_{L^{\frac{q}{q-1}}(V)}\right\}}\cdot{\Lambda_0}^{\max\{p,q\}-1},
\end{eqnarray}
where
$$
\Lambda_0=\min\left\{\frac{l_0}{2}\min\left\{h_0^{\frac{1}{p}}\mu_0^{\frac{1}{p}},h_0^{\frac{1}{q}}\mu_0^{\frac{1}{q}}\right\},1\right\},
$$
system (\ref{eq1}) has one nontrivial solution $(u_{\star,1},v_{\star, 1})$ of  positive energy. Furthermore, if the following condition holds:\\
 $(C_4)$ \;  there exists $l_2>0$, $\beta_3>1$ and $K_3(x):V \to \R$ such that $K_3(x_4)>0$, and $F(x_4,s,s)\geq K_3(x_4)|t|^{\beta_3}$ for  some $x_4\in \{x\in V \big|e_1(x)+e_2(x)>0\}$
 with $\mu(x_4)>0$ and all $s\in \R$ with $0<s< l_2$},\\
 then  system (\ref{eq1}) has another nontrivial solution $(u_{\star,2},v_{\star,2})$ of negative energy  for each  $\lambda\in (0,\lambda_0)$.

\vskip2mm
 \par
 By using similar proofs, we can also obtain some results similar to Theorem 1.1 and Theorem 1.3 to the following equation on locally finite graph $(V,E)$:
 \begin{eqnarray}
\label{ppp1}
   -\Delta_p u+h(x)|u|^{p-2}u=F_u(x,u)+\epsilon e(x),\;\;\;\;\hfill x\in V.
\end{eqnarray}
 \vskip0mm
 \noindent
{\bf Theorem 1.4.} {\it Assume that the following conditions hold:\\
$(h_1)$ \; there exists a constant $h_0>0$ such that $h(x)\geq h_0$ for all $x\in V$;\\
$(h_2)$ \; $h(x)\rightarrow \infty$ as dist$(x,x_0)\rightarrow \infty$ for some fixed $x_0$;\\
$(f_0)$\; $F(x,s)$ is continuously differentiable in $s\in \R$ for all $x\in V$,  and there exists a  function $a\in C(\R^+,\R^+)$ and a function $b:V\to \R^+$ with $b\in L^1(V)$ such that
$$
|F_s(x,s)|\le a(|s|) b(x), \ \  |F(x,s)|\le a(|s|) b(x),
$$
for all $x\in V$ and all $s\in \R$;\\
$(f_1)$ \; $F(x,0)=0$, and there exists $f_1,g_1:V\to \R^+$ with $f_1\in L^\infty(V)$ and $g_1\in L^{\frac{p}{p-1}}(V)$
satisfying $ \|f_1\|_{\infty}<h_{01}$ such that
$$
|F_s(x,s)|\le f_1(x)|s|^{p -1}+g_1(x)\ \ \mbox{for all }x\in V \ \mbox{and all} \ \ s\in \R;
$$
$(f_2)$ there exists $\beta_1>1$ and $K_1:V\to \R$ such that $K_1(x_1)>0$ for some $x_1 \in V$ with $e(x_1)>0$ and
$$
 F(x,s)\ge -K_1(x)|s|^{\beta_1} \ \mbox{for all } x\in V \mbox{ and all  } s\in \R.
$$
 Then for each $\epsilon\in (0,+\infty)$, system (\ref{ppp1})
 has at least  one nontrivial solution.
}

 \vskip2mm
 \noindent
 {\bf Theorem 1.5.} {\it  Assume $(h_1)$, $(f_0)$ and the following conditions hold: \\
$(h_2')$ \; for any given constant $B>0$, $\sum_{x\in A} \mu(x)<\infty$, where $A=\{x\in V|h(x)\le B\}$;\\
$(c_1)$ \; $F(x,0)=0$ for all $x\in V$, and there exists a constant $l_0>0$ such that
$$
|f(x,s)|\le \frac{h_{0}}{p+1}|s|^{p -1}
$$
for all $x\in V$  and all $s\in \R$ with $|s|< l_0$; \\
$(c_2)$ \;  there exists $l_1>0$ such that $F(x_2,s)\geq M s^p\;\mbox{for some } x_2\in V$ with $e(x_2)>0$ and all $s\in \R$ with $s> l_1$, where
$$
M>\frac{D_1+\mu(x_2)h_1(x_2)}{p\mu(x_2)},
$$
 $D_1=\left(\frac{deg(x_2)}{2}\right)^{\frac{p}{2}}\sum_{x\thicksim x_2\mbox{ and }x=x_2}\left(\frac{1}{\mu(x)}\right)^{\frac{p}{2}-1}$;\\
$(c_3)$ \; there exists a constant $\nu>p$ and $0\le A<h_{0}\left(\dfrac{\nu}{p}-1\right)$ such that
$$
\nu F(x,s)-F_s(x,s)s\le A |s|^p \ \ \mbox{for all}\;x \in V.
$$
Then for each $\epsilon$ satisfying
\begin{eqnarray}
\label{pp2}
0<\epsilon< \epsilon_0=\frac{\left(\min\left\{l_0(h_{0}\mu_0)^{\frac{1}{p}}, 1\right\}\right)^{p-1}}{(p+1)h_{0}^{-\frac{1}{p}}\|e\|_{L^{\frac{p}{p-1}}(V)}},
\end{eqnarray}
equation (\ref{ppp1}) has one nontrivial  solution of  positive energy. Furthermore, if the following condition holds:\\
 $(c_4)$ \;  there exists $l_2>0$ and $\beta_3>1$ such that $F(x_3,s)\geq K_3(x_3)|s|^{\beta_3}$ for  some $x_3\in \{x\in V\big|e(x)>0\}$
 and all $s\in \R$ with $0<s< l_2$,\\
 then  system (\ref{ppp1}) has another nontrivial solution  of negative energy  for each  $\epsilon\in (0,\epsilon_0)$.}

\vskip0mm
 \noindent
{\bf Remark 1.1.} In Theorem 1.3, the condition $(C_2)$ is interesting, which implies that the inequality $F(x,s,t)\geq M(s^p+t^q)$ holds only for a
point $x_3$ rather than all $x\in V$ and only for a ray $s=t$ starting at the point $(l_1,l_1)$ in the plane $\R^2$ rather than for all $(s,t)\in \R^2$ with $|(s,t)|>l_1$ (see Fig.1), which is usually assumed in investigating the existence of solutions for the elliptic partial differential system with the nonlinear term satisfying the super-quadratic conditions (for example, see \cite{Liu-Zhang}).
\begin{figure}[h]
 \centering
 \subfigure[Fig.1:  $F(x,s,t)\geq M(s^p+t^q)$ holds only for a ray $s=t$ starting at the point $(l_1,l_1)$ in the plane $\R^2$. ]{
  \includegraphics[scale=0.3]{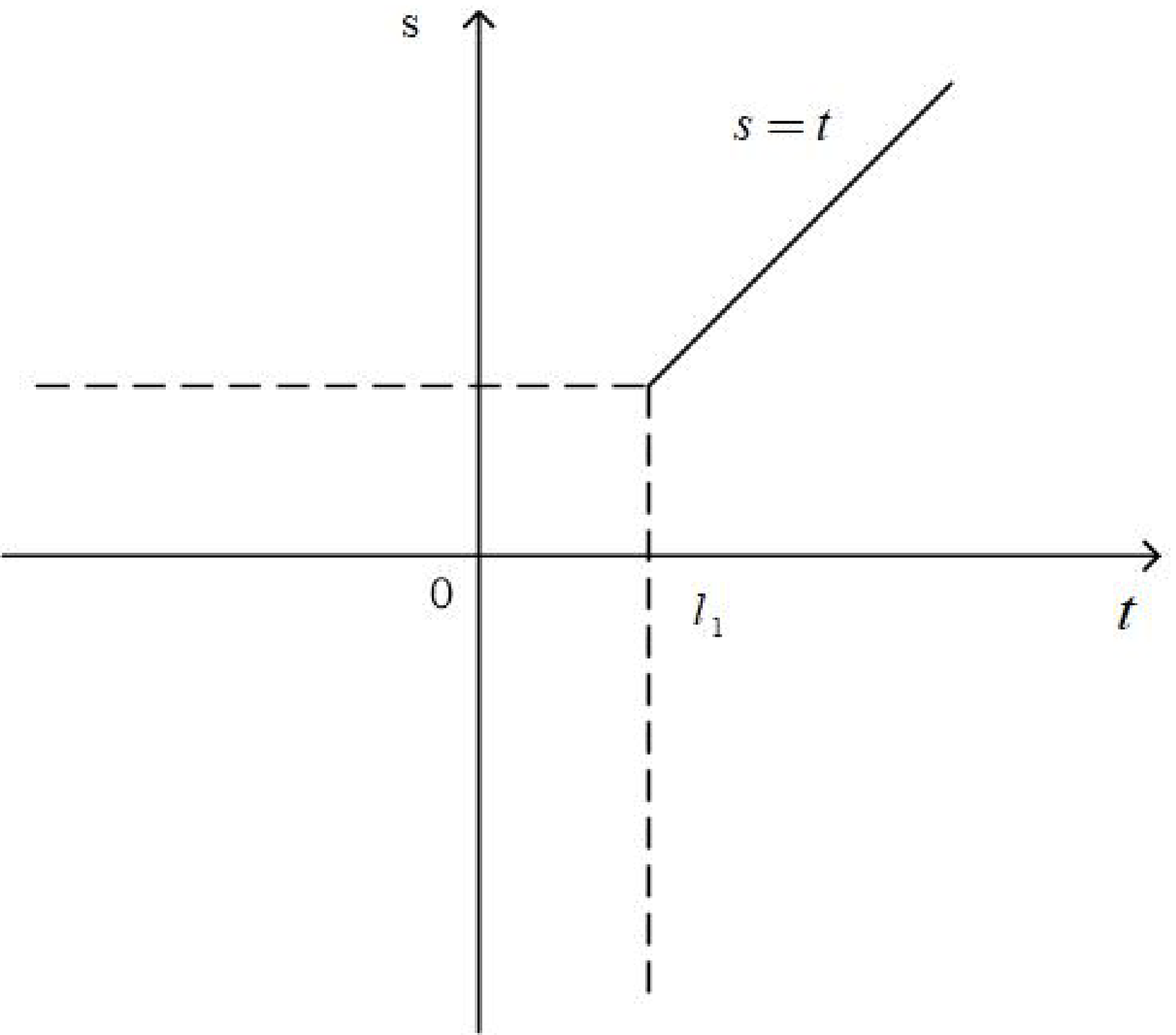}}
\end{figure}
\vskip2mm
 \noindent
{\bf Remark 1.2.} Theorem 1.4 is different from Theorem 1.4 in \cite{Grigor'yan 2017}, where they consider (\ref{ppp1}) with $p=2$ and they assume that $f(x,s):=F_s(x,s)$ satisfies the (AR)-condition $(S_3)$ and $(S_1')$. It is easy to see that  $(f_1)$ in Theorem 1.4 is weaker than $(S_1')$ even if $p=2$ and we do not need the (AR)-condition in Theorem 1.4. Theorem 1.5 is also different from Theorem 1.3 in \cite{Grigor'yan 2017} even if $p=2$.
It is easy to see that $(h_2)$ is weaker than $(K_2)$. Moreover, $(h_2')$ is weaker than $(h_2)$.
In fact, by $(h_2)$ we have for any positive constant $B$ when $h(x)<B$, there exists positive constant $B_1$ such that $dist(x,x_0)<B_1$. So, $A=\{x\in V\big|h(x)<B\}$ is a finite set.
Moreover, $(c_2)$ together with $(c_3)$ is weaker than $(S_3)$. There exists examples satisfying Theorem 1.5 but not satisfying Theorem 1.3 in \cite{Grigor'yan 2017}, for example, let
$$
F(x,s)=M\ln(1+s^2)|s|^3,
$$
where $M$ is defined as Theorem 1.5.

\vskip2mm
\section{Sobolev embedding}
  \setcounter{equation}{0}

Let $W^{1,s}(V)$ be the completion of $C_c(V)$ under the norm
$$
\|u\|_{W^{1,s}(V)}=\left(\int_V|\nabla u(x)|^s+|u(x)|^sd\mu\right)^\frac{1}{s},
$$
where $s>1$ and $W^{1,s}(V)$ is reflexive (see \cite{Meng}).
Let $h(x)\ge h_0>0$. Define the space
$$
W_h^{1,s}(V)=\left\{u\in W^{1,s}(V) \Big|\int_V h(x)|u(x)|^s d\mu<\infty\right\}
$$
endowed with the norm
$$
 \|u\|_{W_h^{1,s}(V)}=\left(\int_V|\nabla u(x)|^s+h(x)|u(x)|^sd\mu\right)^\frac{1}{s}.
$$

  \vskip2mm
  \noindent
{\bf Lemma 2.1.} {\it If $\mu(x)\geq \mu_0>0$ and $h$ satisfies $(H_1)$, then  $W_h^{1,s}(V)$ is continuously embedded into $L^r(V)$ for all $1<s\le r\le \infty$, and the following inequalities hold:
\begin{eqnarray}\label{bb1}
 \|u\|_{\infty}\le \frac{1}{h_0^{1/s}\mu_0^{1/s}} \|u\|_{W_h^{1,s}(V)}
\end{eqnarray}
and
\begin{eqnarray} \label{bb2}
\|u\|_{L^r(V)}\le \mu_0^{\frac{s-r}{sr}}h_0^{-\frac{1}{s}}\|u\|_{W_h^{1,s}(V)}\ \ \mbox{for all } s\le r< \infty.
\end{eqnarray}
Furthermore, if $(H_2)$ also holds, then  $W_h^{1,s}(V)$ is compactly embedded into $L^r(V)$ for all $1<s\le r\le \infty$. }
  \noindent
{\bf Proof.}  For any $u\in W_h^{1,s}(V)$,  we claim that
\begin{eqnarray}\label{aaa5}
\sum_{x\in V} |u(x)|^s \ge \|u\|_\infty^s .
\end{eqnarray}
In fact, assume that
\begin{eqnarray*}
\left(\sum_{x\in V} |u(x)|^s\right)^{1/s} < \|u\|_\infty.
 \end{eqnarray*}
 Then there exists a $\varepsilon>0$ such that
\begin{eqnarray}\label{aaa4}
 \sum_{x\in V} |u(x)|^s < (\|u\|_\infty-\varepsilon)^s.
 \end{eqnarray}
  Note that $\|u\|_\infty=\sup_{x\in V}|u(x)|$. Then by the definition of superum, there exists a $x_*\in V$ such that $|u(x_*)|>\|u\|_\infty-\varepsilon$. Then $|u(x_*)|^s>(\|u\|_\infty-\varepsilon)^s$, which together with (\ref{aaa4}) implies that
  \begin{eqnarray*}
 |u(x_*)|^s>(\|u\|_\infty-\varepsilon)^s >\sum_{x\in V} |u(x)|^s \ge |u(x_*)|^s,
 \end{eqnarray*}
 a contradiction.
 \par
  For any $u\in W_h^{1,s}(V)$,  we have
\begin{eqnarray}\label{aa5}
\|u\|_{W_h^{1,s}(V)}^s\geq\int_V h(x)|u(x)|^s d\mu\geq h_0\int_V |u(x)|^s d\mu \ \ \mbox{for all } s>1,
\end{eqnarray}
 and by (\ref{aaa5}), we have
\begin{eqnarray}\label{aa6}
\|u\|_{W_h^{1,s}(V)}^s\geq\int_V h(x)|u(x)|^s d\mu = \sum_{x\in V}\mu(x)h(x)|u(x)|^s \geq h_0\mu_0 \|u\|_{\infty}^s \ \ \mbox{for all } s>1,
\end{eqnarray}
which implies that
$$
 \|u\|_{\infty}\le \frac{1}{h_0^{1/s}\mu_0^{1/s}} \|u\|_{W_h^{1,s}(V)} \ \ \mbox{for all } s>1.
$$
When $s<r<\infty$, it follows from (\ref{aa5}) and (\ref{aa6}) that
$$
\int_V |u(x)|^r d\mu\le\|u\|_{\infty}^{r-s} \int_V|u(x)|^{s} d\mu\le \mu_0^{\frac{s-r}{s}}h_0^{-\frac{r}{s}}\|u\|_{W_h^{1,s}(V)}^r.
$$
So,
$$
\|u\|_{L^r(V)}\le \mu_0^{\frac{s-r}{sr}}h_0^{-\frac{1}{s}}\|u\|_{W_h^{1,s}(V)} \ \ \mbox{for all } s\le r< \infty.
$$

\par
Suppose that $\{u_k\}$ is a bounded sequence in $W_h^{1,s}(V)$. Note that $W_h^{1,s}(V)$ is reflexive. Then there exists a subsequence, still denoted by $\{u_k\}$, such that $u_k \rightharpoonup u$ weakly in $W_h^{1,s}(V)$ for some $u\in W_h^{1,s}(V)$. In particular,
$$
\lim_{k\rightarrow\infty}\int_V u_k\varphi d\mu=\int_V u\varphi d\mu,\;\forall\varphi\in C_c(V),
$$
which implies that
\begin{eqnarray} \label{aaa1}
\lim_{k \rightarrow \infty}u_k(x)=u(x) \ \ \mbox{for any fixed } x\in V,
\end{eqnarray}
if we choose $\varphi\in C_c(V)$ defined by
$$
\varphi(y)=\begin{cases}
1,& y= x,\\
0,& y\not=x.
\end{cases}
$$
 We now prove $u_k \rightarrow u$ in $L^r(V)$ for all $s \le r\le \infty$, if $(H_2)$ holds. Since $\{u_k\}$ is bounded in $W_h^{1,s}(V)$ and $u \in W_h^{1,s}(V)$, by the definition of norm $\|\cdot\|_{W_h^{1,s}(V)}$, there exists a constant $c_0>0$ such that
$$
\int_V h|u_k-u|^s d\mu\le c_0.
$$
For any given $\epsilon>0$, in view of $(H_2)$, there exists a constant $R(\epsilon)>0$ such that
$$
\frac{1}{h(x)}<\epsilon\;\; \mbox{as dist}(x,x_0)\geq R(\epsilon).
$$
Hence,
\begin{eqnarray}\label{aaa2}
\int_{\mbox{dist}(x,x_0)\geq R(\epsilon)} |u_k-u|^s d\mu&=& \int_{\mbox{dist}(x,x_0)\geq R(\epsilon)} \frac{1}{h} h|u_k-u|^s d\mu\nonumber\\
                                     &\le& c_0\epsilon.
\end{eqnarray}
Note that $\{x|\mbox{dist}(x,x_0)\le R(\epsilon)\}$ is a finite set. Then (\ref{aaa1}) implies that
\begin{eqnarray}\label{aaa3}
\lim_{k\rightarrow \infty}\int_{\mbox{dist}(x,x_0)\le R(\epsilon)} |u_k-u|^s d\mu =0.
\end{eqnarray}
So, by the arbitrary of $\epsilon$, (\ref{aaa2}) and (\ref{aaa3}) imply that
\begin{eqnarray}\label{aaa6}
\lim_{k\rightarrow \infty}\int_V |u_k-u|^s d\mu =0.
\end{eqnarray}
Then by (\ref{aaa5}) and (\ref{aaa6}), we have
\begin{eqnarray}\label{bbb1}
         \|u_k-u\|_{\infty}^s
 & \le & \sum_{x\in V}|u_k(x)-u(x)|^s\nonumber\\
 &  =  & \sum_{x\in V}\frac{1}{\mu(x)}\mu(x)|u_k(x)-u(x)|^s\nonumber\\
 & \le & \frac{1}{\mu_0}\sum_{x\in V}\mu(x)|u_k(x)-u(x)|^s\nonumber\\
 &  =  & \frac{1}{\mu_0} \int_V |u_k-u|^s d\mu \to 0, \ \ \mbox{as }k\to \infty,
\end{eqnarray}
and when $s<r<\infty$, we have
\begin{eqnarray}\label{aaa7}
\int_V |u_k-u|^r d\mu\le\|u_k-u\|_{\infty}^{r-s}\int_V |u_k-u|^s d\mu\to 0, \ \ \mbox{as }k\to \infty.
\end{eqnarray}
Hence, (\ref{aaa6}), (\ref{bbb1})  and (\ref{aaa7}) imply that $u_k\rightarrow u$ in $L^r(V)$ for all $s\le r\le\infty$.
\qed

\vskip2mm
\noindent
{\bf Lemma 2.2.} {\it If $\mu(x)\geq \mu_0>0$ and $h$ satisfies $(H_1)$ and $(H_2')$, then $W_h^{1,s}(V)$ is compactly embedded into $L^r(V)$ for all $1<s\le r\le \infty.$}
\vskip2mm
\noindent
{\bf Proof.}
 Suppose that $\{u_k\}$ is a bounded sequence in $W_h^{1,s}(V)$ and there exists a positive constant $C_0$ such that
\begin{eqnarray}\label{aaaa2}
 \|u_k\|_{W_h^{1,s}(V)}\le C_0.
\end{eqnarray}
 Since $\|u_k\|_{L^s(V)}^s\le \frac{1}{h_0}\int_V h(x)|u_k|^s d\mu \le \frac{1}{h_0}\|u_k\|_{W_h^{1,s}(V)}^s$,
 we also have $\{\|u_k\|_{L^s(V)}\}$ is bounded  in $\R$.
 Noting that $W_h^{1,s}(V)$ is reflexive, we have, up to a subsequence, $u_k \rightharpoonup u$ weakly in $W_h^{1,s}(V)$ for some $u\in W_h^{1,s}(V)$ and $\delta_k=\|u_k\|_{L^s(V)} \rightarrow \delta$ for some $\delta\in \R$ as $k\to \infty$. Similar to the argument of (\ref{aaa1}), we have $\lim_{k\to \infty} u_k(x)=u(x)$ for all $x\in V$. Then for any bounded domain $\Omega \subset V$, we have
$$
 \int_\Omega |u_k|^s d\mu\rightarrow \int_\Omega |u|^s d\mu  \ \mbox{and} \int_\Omega |u_k|^s d\mu\le \int_V |u_k|^s d\mu \rightarrow \delta^s, \ \mbox{ as } k\to\infty.
$$
Then
\begin{eqnarray}\label{aaaa3}
\delta^s\geq \|u\|_{L^s(\Omega)}^s.
\end{eqnarray}
For any given constant $B>0$,  define $\Omega=\{x\in V|\mbox{dist}(x,x_0)\le B,h(x)\le B\}$ for some fixed $x_0\in V$. Let $A(\Omega)=\{x\in V/\Omega, h(x)\le B\}$. Then $A=\Omega\cup A(\Omega)$, where $A=\{x\in V|h(x)\le B\}$. By $(H_2')$, we have $\sum_{x\in A} \mu(x)<\infty$ and then by the definition of convergent series, for any sufficient small $\epsilon>0$, there exists a sufficient large $B>\frac{1}{\epsilon}$ such that
\begin{eqnarray}\label{aaaa1}
\sum_{x\in A(\Omega)} \mu(x)<\epsilon.
\end{eqnarray}
Moreover, since $h$ satisfies $(H_1)$, by Lemma 2.1 we know that $W_h^{1,s}(V)$ is continuously embedded into $L^r(V)$, $s\le r\le \infty$. So, by (\ref{bb2}), (\ref{aaaa2}) and (\ref{aaaa1}), we have
\begin{eqnarray*}
\int_{A(\Omega)} |u_k|^s d\mu & =  &  \int_{A(\Omega)}1\cdot |u_k|^s d\mu\\
                              &\le& \left(\int_{A(\Omega)}|u_k|^{2s}d\mu\right)^{\frac{1}{2}}\left(\sum_{x\in A(\Omega)}\mu(x)\right)^{\frac{1}{2}}\\
                              &\le& \mu_0^{-\frac{1}{2}}h_0^{-1}\|u_k\|_{W_h^{1,s}(V)}^s\left(\sum_{x\in A(\Omega)}\mu(x)\right)^{\frac{1}{2}}\\
                              &\le& \mu_0^{-\frac{1}{2}}h_0^{-1}C_0^s\epsilon\ \ \ \mbox{for all } k\in \mathbb N.
\end{eqnarray*}
 Define $B(\Omega)=\{x\in V/\Omega,h(x)>B\}$.  Then
$$
\int_{B(\Omega)} |u_k|^s d\mu\le\int_{B(\Omega)} \frac{h(x)}{B}|u_k|^sd\mu\le \frac{1}{B}\|u_k\|_{W_h^{1,s}(V)}^s\le \frac{C_0^s}{B}\le C_0^s \epsilon \ \ \mbox{for all } k\in \mathbb N.
$$
Then
\begin{eqnarray}\label{aaaa4}
\int_{V/\Omega} |u_k|^s d\mu=\int_{B(\Omega)} |u_k|^s d\mu+\int_{A(\Omega)} |u_k|^s d\mu <\left(\mu_0^{-\frac{1}{2}}h_0^{-1}+1\right)C_0^s\epsilon\ \ \ \mbox{for all } k\in \mathbb N.
\end{eqnarray}
Similarly, we also have
\begin{eqnarray}\label{aaaa5}
\int_{V/\Omega} |u|^s d\mu=\int_{B(\Omega)} |u|^s d\mu+\int_{A(\Omega)} |u|^s d\mu <\epsilon\left(\mu_0^{-\frac{1}{2}}h_0^{-1}+1\right)\|u\|_{W_h^{1,s}(V)}^s.
\end{eqnarray}
Let $C_1=\max\left\{C_0,\|u\|_{W_h^{1,s}(V)}\right\}$. So, by (\ref{aaaa3}) and  (\ref{aaaa5}), we have
\begin{eqnarray*}
 \|u\|_{L^s(V)}^s
&   =   & \|u\|_{L^s(\Omega)}^s+\|u\|_{L^s(V/\Omega)}^s\\
& \le  & \delta^s+\left(\mu_0^{-\frac{1}{2}}h_0^{-1}+1\right)C_1^s\epsilon.
\end{eqnarray*}
On the other hand,
\begin{eqnarray*}
         \|u\|_{L^s(V)}^s
&   =   & \|u\|_{L^s(\Omega)}^s+\|u\|_{L^s(V/\Omega)}^s\\
& \ge  & \lim_{k\rightarrow \infty}\|u_k\|_{L^s(\Omega)}^s\\
&  =    &  \lim_{k\rightarrow \infty}\|u_k\|_{L^s(V)}^s-\lim_{k\rightarrow \infty}\|u_k\|_{L^s(V/\Omega)}^s\\
& \ge  & \delta^s-\left(\mu_0^{-\frac{1}{2}}h_0^{-1}+1\right)C_1^s\epsilon.
\end{eqnarray*}
Hence, by the arbitrary of $\epsilon$, we obtain that $\delta^s= \|u\|_{L^s(V)}^s$. Thus we have proved that
$\|u_k\|_{L^s(V)} \to \|u\|_{L^s(V)}$ as $k\to \infty$. By the uniform convexity of $L^s(V)$ (Lemma 2.2 in (\cite{Meng})) and that $u_k \rightharpoonup u$ weakly in $W_h^{1,s}(V)$, it follows from the Kadec-Klee property that $\|u_k-u\|_{L^s(V)}\to 0$ as $k\to \infty$. Then similar to the argument of (\ref{bbb1}) and (\ref{aaa7}), we have $\|u_k-u\|_{\infty}\to 0$ and $\|u_k-u\|_{L^r(V)}\to 0$ for all $s<r<\infty$.\qed

\vskip2mm
\noindent
{\bf Remark 2.1.} Lemma 2.1 generalizes Lemma 2.2  in \cite{Grigor'yan 2017} and  Lemma 2.6 in \cite{Meng}, and Lemma 2.2 generalizes  Lemma 3 in \cite{Chang2022}. To be precise, when $s=2$,
Lemma 2.1 and Lemma 2.2 reduce to Lemma 2.2 in \cite{Grigor'yan 2017}
and  Lemma 3 in \cite{Chang2022}, respectively.
In Lemma 3 in \cite{Chang2022}, the potential $h(x)$ is allowed to be sigh-changing,
which satisfies $(H_1')$: $\inf_{x\in V} h(x)\ge h_0$ for some $h_0\in (-1,0)$.
One can prove that Lemma 2.2 still holds under $(H_1')$ and $(H_2')$. Moreover, if $h(x)=\lambda a(x)+1$, where $a:V\to\R$ with $a(x)\ge 0$ for all $x\in V$, then Lemma 2.1 reduces to Lemma 2.6 in \cite{Meng}.
The proofs of Lemma 2.1 and Lemma 2.2 are based on those  in \cite{Grigor'yan 2017}, \cite{Chang2022} and \cite{Meng} and we make some appropriate modifications.

\vskip2mm
\par
Assume that $\varphi \in C^{1}(X, \mathbb{R})$. An sequence $\{u_n\}$ is called as the Palais-Smale sequence of $\varphi$ if $\varphi(u_n)$ is bounded for all $n\in \mathbb N$ and  $\varphi'(u_n)\to 0$ as $n\to \infty$. If any Palais-Smale sequence $\{u_n\}$ of $\varphi$ has a convergent subsequence, we call that $\varphi$ satisfies the Palais-Smale condition ((PS)-condition for short).

\vskip2mm
\noindent
{\bf Lemma 2.3.} (Ekeland's variational principle \cite{Mawhin J 1989}) {\it Let $M$ be a  complete metric space with metric $d$, and $\varphi:M\rightarrow \R$ be a lower semicontinuous function, bounded from below and not identical to $+\infty$. Let $\varepsilon >0$ be given and $U \in M$ such that
$$
\varphi(U)\le \inf_{M} \varphi+\varepsilon.
$$
Then there exists $V \in M$ such that
$$
\varphi(V)\le \varphi(U),\;d(U,V)\le1,
$$
and for each $W\in M$, one has
$$
\varphi(V)\le \varphi(W)+\varepsilon d(V,W).
$$}

\par
By the Ekeland's variational principle, it is easy to obtain the following corollary.
\vskip2mm
\noindent
{\bf Lemma 2.4.} (\cite{Mawhin J 1989}) {\it Suppose that $X$ is a Banach space, $M\subset X$ is closed, $\varphi\in C^1(X,\R)$ is bounded from below on $M$ and satisfies the (PS)-condition. Then $\varphi$ attains its infimum on $M$.}
\vskip2mm
\noindent
{\bf Lemma 2.5.} (Mountain pass theorem \cite{Rabinowitz 1986}) {\it Let $X$ be a real Banach space and $\varphi \in C^{1}(X,\R)$, $\varphi(0)=0$ satisfy (PS)-condition. Suppose that $\varphi $ satisfies the following conditions:\\
(i) there exists a constant $ \rho>0$ and $\alpha>0$ such that $ \varphi|_{\partial B_{\rho}(0)}\ge \alpha $, where $B_\rho=\{w\in X:\|w\|_X<\rho\}$;\\
(ii) there exists $ w\in X\backslash \bar B_{\rho} (0)$ such that $ \varphi(w)\le 0 $.\\
Then $\varphi$ has a critical value $c_*\ge \alpha$ with
 $$
 c_*:=\inf_{\gamma\in\Gamma}\max_{t\in[0,1]}\varphi(\gamma(t)),
 $$
where
 $$
 \Gamma:=\{\gamma\in C([0,1],X]):\gamma(0)=0,\gamma(1)=w\}.
 $$
 }

\vskip2mm
\section{Proofs for the sub-$(p,q)$-linear case}
  \setcounter{equation}{0}
  \par
Define the space $W:=W_h^{1,p}(V)\times W_h^{1,q}(V)$ with the norm
$$
\|(u,v)\|_W=\|u\|_{W_h^{1,p}(V)}+\|v\|_{W_h^{1,q}(V)}.
$$
Then $W$ is a Banach space. Consider the functional $\varphi:W\to\R$ defined as
\begin{eqnarray}
\label{eq2}
\notag \varphi_\lambda(u,v)&=&\frac{1}{p}\int_V(|\nabla u|^p+h_1|u|^p)d\mu+\frac{1}{q}\int_V(|\nabla v|^q+h_2|v|^q)d\mu\\
                   &&-\int_V F(x,u,v)d\mu-\lambda_1\int_Ve_1u d\mu-\lambda_2\int_Ve_2vd\mu.
\end{eqnarray}
Then $\varphi_\lambda(u,v)\in C^1(W,\R)$, and
\begin{eqnarray}
\label{eq3}
\notag \langle\varphi_\lambda'(u,v),(\phi_1,\phi_2)\rangle&=&\int_V\left[|\nabla u|^{p-2}\Gamma(u,\phi_1)+h_1|u|^{p-2}u\phi_1-F_u(x,u,v)\phi_1-\lambda_1 e_1\phi_1\right]d\mu\\
&&+\int_V\left[|\nabla v|^{q-2}\Gamma(v,\phi_2)+h_2|v|^{q-2}v\phi_2-F_v(x,u,v)\phi_2-\lambda_2 e_2\phi_2\right]d\mu
\end{eqnarray}
for all $(\phi_1,\phi_2)\in W$ (see Appendix A.2).
\vskip2mm
\noindent
{\bf Definition 3.1. }{\it  $(u,v)\in W$ is called as a weak solution of system (\ref{eq1}) if
\begin{eqnarray}
 \label{pp10}&   &       \int_V\left[|\nabla u|^{p-2}\Gamma(u,\phi_1)+h_1|u|^{p-2}u\phi_1\right]d\mu
         =   \int_V[F_u(x,u,v)\phi_1+\lambda_1 e_1\phi_1]d\mu,\\
\label{pp11} &  &    \int_V\left[|\nabla v|^{q-2}\Gamma(v,\phi_2)+h_2|v|^{q-2}v\phi_2\right]d\mu
        =   \int_V[F_v(x,u,v)\phi_1+\lambda_2e_2\phi_2]d\mu,
\end{eqnarray}
for all $(\phi_1,\phi_2)\in W$.
}
\par
Obviously, $(u,v)\in W$ is  a weak solution of system (\ref{eq1}) if and only if $(u,v)$ is a critical point of $\varphi$ and similar to the arguments in \cite{Meng}, we have the
following proposition.
\vskip2mm
\noindent
{\bf Proposition 3.1. }{\it If  $(u,v)\in W$ is a weak solution of system (\ref{eq1}), then $(u,v)\in W$ is also a point-wise solution of (\ref{eq1}). }
\vskip2mm
\noindent
{\bf Proof.} For any
fixed $y\in V$, we take a test function $\phi_1:V\to \R$ in (\ref{pp10}) with
\begin{eqnarray*}
\phi_1(x)=\begin{cases}
1,& x=y,\\
0,& x\not=y,
\end{cases}
\end{eqnarray*}
and a test function $\phi_2:V\to \R$ in (\ref{pp11}) with
\begin{eqnarray*}
\phi_2(x)=\begin{cases}
1,& x=y,\\
0,& x\not=y.
\end{cases}
\end{eqnarray*}
Thus, by (\ref{p1}), we have
\begin{eqnarray*}
&   &       -\Delta_pu(y)+h_1(y)|u(y)|^{p-2}u(y)
         =   F_u(y,u(y),v(y))+\lambda_1 e_1(y),\\
 &  &     -\Delta_qv(y)+h_2(y)|v(y)|^{q-2}v(y)
        =   F_v(y,u(y),v(y))+\lambda_2e_2(y).
\end{eqnarray*}
By the arbitrary of $y$, we complete the proof.\qed
\vskip2mm
\noindent
{\bf Lemma 3.1.} {\it Assume that $(H_1)$ and $(F_1)$ hold. Then $\varphi$ is coercive, that is, $\varphi(u,v) \to +\infty$ as $\|(u,v)\|_W\to \infty$.}
\vskip2mm
\noindent
{\bf Proof.} By $(F_1)$ and Lemma 2.1, we have
\begin{eqnarray}
\label{EQ15}
\notag \int_V|F(x,u,v)|d\mu &=&\int_V|F(x,u,v)-F(x,0,0)|d\mu\\
\notag &\le& \int_V|F(x,u,v)-F(x,0,v)|+|F(x,0,v)-F(x,0,0)|d\mu\\
\notag &\le& \int_V \int_0^{|u|}|F_s(x,s,v)|dsd\mu+\int_V \int_0^{|v|}|F_t(x,0,t)|dtd\mu\\
\notag &\le& \int_V\int_0^{|u|}[f_1(x)(|s|^{p-1}+|v|^{\frac{pq-q}{p}})+g_1(x)]dsd\mu+\int_0^{|v|}[ f_2(x)|t|^{q-1}+g_2(x)]dtd\mu\\
\notag &\le& \int_V \left[\frac{|u|^p}{p}f_1(x)+f_1(x)|v|^{\frac{pq-q}{p}}|u|+\frac{|v|^q}{q}f_2(x)+g_1(x)|u|+g_2(x)|v| \right]d\mu\\
\notag &\le& \frac{2\|f_1\|_{\infty}}{p}\int_V|u|^pd\mu+\frac{(p-1)\|f_1\|_{\infty}}{p}\int_V|v|^qd\mu\\
\notag & &    +\frac{\|f_2\|_{\infty}}{q}\int_V|v|^qd\mu+\|g_1\|_{L^{\frac{p}{p-1}}(V)}\|u\|_{L^p(V)}+\|g_2\|_{L^{\frac{q}{q-1}}(V)}\|v\|_{L^q(V)}\\
\notag&\le&\frac{2\|f_1\|_{\infty}}{ph_0}\|u\|_{W_h^{1,p}(V)}^p+\left(\frac{(p-1)\|f_1\|_{\infty}}{ph_0}+\frac{\|f_2\|_{\infty}}{qh_0}\right)\|v\|_{W_h^{1,q}(V)}^q\\
       &&     +\frac{\|g_1\|_{L^{\frac{p}{p-1}}(V)}}{{h_0}^{\frac{1}{p}}}\|u\|_{W_h^{1,p}(V)}+\frac{\|g_2\|_{L^{\frac{q}{q-1}}(V)}}{{h_0}^{\frac{1}{q}}}\|v\|_{W_h^{1,q}(V)}.
\end{eqnarray}
Then, by (\ref{eq2}) and (\ref{EQ15}), we have
\begin{eqnarray*}
\varphi_\lambda(u,v)
& \ge & \left(\frac{1}{p}-\frac{2\|f_1\|_{\infty}}{ph_0}\right)\|u\|_{W_h^{1,p}(V)}^p
         +\left(\frac{1}{q}-\frac{(p-1)\|f_1\|_{\infty}}{ph_0}-\frac{\|f_2\|_{\infty}}{qh_0}\right)\|v\|_{W_h^{1,q}(V)}^q\\
&&       -\frac{1}{{h_0}^{\frac{1}{p}}}\left(\lambda_1\|e_1\|_{L^{\frac{p}{p-1}}(V)}+\|g_1\|_{L^{\frac{p}{p-1}}(V)}\right)\|u\|_{W_h^{1,p}(V)}\\
&&       -\frac{1}{{h_0}^{\frac{1}{q}}}\left(\lambda_2\|e_2\|_{L^{\frac{q}{q-1}}(V)}+\|g_2\|_{L^{\frac{q}{q-1}}(V)}\right)\|v\|_{W_h^{1,q}(V)}.
\end{eqnarray*}
So $\varphi$ is coercive in $W$.\qed

\vskip2mm
\noindent
{\bf Lemma 3.2.} {\it Assume that $(H_1)$ and $(F_1)$  hold. Then $\varphi_{\lambda}$  satisfies the (PS)-condtion.}

\vskip2mm
\noindent
{\bf Proof.} The proof is motivated by \cite{Meng} and \cite{Zhang-Tang}. Assume that $\{(u_k,v_k)\}$ is a Palais-Smale sequence, then $\varphi'_{\lambda}(u_k,v_k)\rightarrow 0$ as $k\rightarrow\infty$ and $\varphi_{\lambda}(u_k,v_k)$ is bounded.
 By Lemma 3.1, we obtain that $\{(u_k,v_k)\}$ is bounded in $W$.
Then $\{u_k\}$ is bounded in $W_h^{1,p}(V)$ and $\{v_k\}$ is bounded in $W_h^{1,q}(V)$. Hence we can find a
subsequence, still denoted by $\{u_k\}$, such that $u_k\rightharpoonup u_{\lambda \star}$ for some $u_{\lambda \star}\in W_h^{1,p}(V)$ as $k\to \infty$, and a
subsequence of $\{v_k\}$,  which has the same subscript as the subsequence of $\{u_k\}$, still denoted by $\{v_k\}$,
such that $v_k\rightharpoonup v_{\lambda \star}$ for some $v_{\lambda \star}\in W_h^{1,q}(V)$ as $k\to \infty$.
By Lemma 2.1, we know that
\begin{eqnarray}
\label{L3.2.1} &  &  u_k\to  u_{\lambda \star} \ \mbox{in } L^p(V),\ \ v_k\to v_{\lambda \star}  \ \mbox{in } L^q(V), \ \ \mbox{as }k\to \infty.
\end{eqnarray}
   Then by (\ref{eq3}), we have
\begin{eqnarray}
\label{aEQ9}
\notag &   & \langle\varphi_\lambda'(u_k,v_k)-\varphi_\lambda'(u_{\lambda \star},v_{\lambda \star}),(u_k-u_{\lambda \star},0)\rangle\\
\notag & = & \int_V\left[|\nabla u_k|^{p-2}\Gamma(u_k,u_k-u_{\lambda \star})+(h_1(x)|u_k|^{p-2}u_k-F_u(x,u_k,v_k))(u_k-u_{\lambda \star})\right]d\mu\\
\notag &   & -\int_V\left[|\nabla u_{\lambda \star}|^{p-2}\Gamma(u_{\lambda \star},u_k-u_{\lambda \star})+(h_1(x)|u_{\lambda \star}|^{p-2}u_{\lambda \star}-F_u(x,u_{\lambda \star},v_{\lambda \star}))(u_k-u_{\lambda \star})\right]d\mu\\
\notag & = & \|u_k\|_{W_h^{1,p}(V)}^p+\|u_{\lambda \star}\|_{W_h^{1,p}(V)}^p-\int_V \left[|\nabla u|^{p-2}\Gamma(u_k, u_{\lambda \star})+h_1(x)|u_k|^{p-2}u_ku_{\lambda \star}\right]d\mu\\
\notag &   & -\int_V \left[|\nabla u|^{p-2}\Gamma(u_{\lambda \star}, u_k)+h_1(x)|u_{\lambda \star}|^{p-2}u_{\lambda \star}u_k \right]d \mu+\int_V[F_u(x,u_{\lambda \star},v_{\lambda \star})-F_u(x,u_k,v_k)](u_k-u_{\lambda \star})d\mu.\\
\end{eqnarray}
By $(F_1)$ and (\ref{L3.2.1}), we have
\begin{eqnarray}
\label{aEQ10}
\notag &&     \int_V[F_u(x,u_{\lambda \star},v_{\lambda \star})-F_u(x,u_k,v_k)](u_k-u_{\lambda \star})d\mu\\
\notag &\le& \int_V|F_u(x,u_{\lambda \star},v_{\lambda \star})-F_u(x,u_k,v_k)||u_k-u_{\lambda \star}|d\mu\\
\notag &\le& \int_V[|F_u(x,u_k,v_k)|+|F_u(x,u_{\lambda \star},v_{\lambda \star})|]|u_k-u_{\lambda \star}|d\mu\\
\notag &\le& \|f_1\|_\infty\int_V(|u_k|^{p-1}+|v_k|^{\frac{pq-q}{p}}+|u_{\lambda \star}|^{p-1}+|v_{\lambda \star}|^{\frac{pq-q}{p}})|u_k-u_{\lambda \star}|d\mu
              +\int_Vg_1(x)|u_k-u_{\lambda \star}|d\mu\\
\notag &\le& \|f_1\|_\infty\left(\|u_k\|_{L^p(V)}^{p-1}+\|v_k\|_{L^q(V)}^{\frac{pq-q}{p}}+\|u_{\lambda \star}\|_{L^p(V)}^{p-1}
              +\|v_{\lambda \star}\|_{L^q(V)}^{\frac{pq-q}{p}}\right)\|u_k-u_{\lambda \star}\|_{L^p(V)}\\
\notag &  &   +\|g_1\|_{L^{\frac{p}{p-1}}(V)}\|u_k-u_{\lambda \star}\|_{L^p(V)}\\
       &\rightarrow& 0.
\end{eqnarray}
Moreover, by (\ref{G1}), we have
\begin{eqnarray}
\label{aEQ11}
\notag &     & \int_V \left[|\nabla u_k|^{p-2}\Gamma(u_k, u_{\lambda \star})+h_1|u_k|^{p-2}u_ku_{\lambda \star}\right]d\mu\\
\notag & \le & \int_V |\nabla u_k|^{p-2}|\nabla u_k ||\nabla  u_{\lambda \star}| d\mu+\int_V (h_1^{\frac{p-1}{p}}|u_k|^{p-2}u_k)(h_1^{\frac{1}{p}}u_{\lambda \star}) d\mu\\
\notag & \le & \|\nabla u_k\|_{L^p(V)}^{p-1}\|\nabla  u_{\lambda \star}\|_{L^p(V)}+\left(\int_V h|u_k|^p d\mu\right)^{\frac{p-1}{p}}\left(\int_Vh|u_{\lambda \star}|^p d\mu\right)^{\frac{1}{p}}\\
       & \le & \|u_k\|_{W_h^{1,p}(V)}^{p-1}\|u_{\lambda \star}\|_{W_h^{1,p}(V)}.
\end {eqnarray}

Similarly, we also have $\int_V \left[|\nabla u_{\lambda \star}|^{p-2}\Gamma(u_{\lambda \star}, u_k)+h_1(x)|u_{\lambda \star}|^{p-2}u_{\lambda \star}u_k \right]d \mu\le \|u_{\lambda \star}\|_{W_h^{1,p}(V)}^{p-1}\|u_k\|_{W^{{m_1},p}(V)}$.
So, by (\ref{aEQ9}), (\ref{aEQ10}) and (\ref{aEQ11}), we have
\begin{eqnarray*}
&&    \langle\varphi_\lambda'(u_k,v_k)-\varphi_\lambda'(u_{\lambda \star},v_{\lambda \star}),(u_k-u_{\lambda \star},0)\rangle\\
&\ge&\|u_k\|_{W_h^{1,p}(V)}^p+\|u_{\lambda \star}\|_{W_h^{1,p}(V)}^p-\|u_k\|_{W^{{m_1},p}(V)}^{p-1}\|u_{\lambda \star}\|_{W_h^{1,p}(V)}-\|u_{\lambda \star}\|_{W_h^{1,p}(V)}^{p-1}\|u_k\|_{W_h^{1,p}(V)}+o_k(1)\\
&=&   (\|u_k\|^{p-1}_{W_h^{1,p}(V)}-\|u_{\lambda \star}\|^{p-1}_{W_h^{1,p}(V)})(\|u_k\|_{W_h^{1,p}(V)}-\|u_{\lambda \star}\|_{W_h^{1,p}(V)})+o_k(1).
\end{eqnarray*}
Hence,  $\|u_k\|_{W_h^{1,p}(V)}\to\|u_{\lambda \star}\|_{W_h^{1,p}(V)}$ as $k\to \infty$. Then it follows from the uniformly convexity of $W_h^{1,p}(V)$ (see Appendix A.1) and the Kadec-Klee property that
$$
u_k \rightarrow u_{\lambda \star} \mbox{ strongly in}\;W_h^{1,p}(V)\; \mbox{as} \;k\rightarrow \infty.
$$
Similarly, we can also prove
$$
v_k\rightarrow v_{\lambda \star} \; \mbox{strongly in} \;W_h^{1,q}(V)\;\mbox{as} \; k\rightarrow \infty.
$$
Therefore,
$$
(u_k,v_k)\rightarrow (u_{\lambda \star},v_{\lambda \star}) \; \mbox{strongly   in}\; W \;\mbox{as}\; k\rightarrow \infty.
$$
\qed

\vskip2mm
\noindent
{\bf Proof of Theorem 1.1.}  By Lemma 3.1 and the continuity of $\varphi_{\lambda}$, we know that $\varphi_{\lambda}$ is bounded from below.
 Then by  Lemma 3.2 and Lemma 2.4, we obtain that $\varphi_{\lambda}$ attains its infimum on $W$.
  Hence, there exists a $(u_{\lambda \star},v_{\lambda \star})\in W$ such that $\varphi(u_{\lambda \star},v_{\lambda \star})=\inf_{(u,v)\in W}\varphi(u,v)$.
 \par
 Next, we prove $(u_{\lambda \star},v_{\lambda \star})\not=(0,0)$.
 Assume that $(u_{\lambda \star},v_{\lambda \star})=(0,0)$. Then
$
\varphi(0,0)=0=\inf_{(u,v)\in W}\varphi(u,v).
$
 Let
 $$
u_*(x)=\begin{cases}
1,& x= x_1,\\
0,& x\not=x_1,
\end{cases}
$$
where $x_1 \in V\mbox{ with }e_1(x_1)>0$. If (i) of $(F_2)$ holds, then
\begin{eqnarray}\label{bb3}
          \inf_{(u,v)\in W}\varphi_\lambda(u,v)
& \le  &  \inf_{\theta\in (0,+\infty)}\varphi_\lambda(\theta u_*,0)\nonumber\\
&  =   &  \inf_{\theta\in (0,+\infty)} \left(\frac{1}{p}\theta^p\|u_*\|_{W_h^{1,p}(V)}^p-\int_V F(x,\theta u_*,0)d\mu-\lambda_1\theta\int_Ve_1u_* d\mu\right)\nonumber\\
& \le &  \inf_{\theta\in (0,+\infty)} \left(\frac{1}{p}\theta^p\|u_*\|_{W_h^{1,p}(V)}^p+\int_V K_1(x)|\theta u_*|^{\beta_1}d\mu-\lambda_1\theta\int_Ve_1u_* d\mu\right)\nonumber\\
&  =  &  \inf_{\theta\in (0,+\infty)} \left(\frac{1}{p}\theta^p\|u_*\|_{W_h^{1,p}(V)}^p+\mu(x_1)\theta ^{\beta_1}K_1(x_1)-\lambda_1\theta\mu(x_1)e_1(x_1)\right).
\end{eqnarray}
Note that $\beta_1>1$, $p>1$, $\mu(x_1)>0$, $e_1(x_1)>0$, $K_1(x_1)>0$ and $\lambda_1>0$. Then for each $\lambda_1>0$, there exists sufficiently
 small $\theta>0$ such that $ \inf_{(u,v)\in W}\varphi(u,v)<0$, which is a contradiction. Similarly,
 if (ii) of $(F_2)$ holds, we also can obtain the same contradiction.
 \par
Moreover, if $(u_{\lambda \star},v_{\lambda \star})=(u_{\lambda \star},0)$, then by (\ref{pp10}), we have
$$
\int_V(|\nabla u_{\lambda \star}|^p+h_1|u_{\lambda \star}|^p)d\mu=\int_V F_u(x,u_{\lambda \star},0)u_{\lambda \star}d\mu+\lambda_1\int_Ve_1u_{\lambda \star} d\mu.
$$
Hence, combining with $(F_1)$, we have
$$
\|u_{\lambda \star}\|_{W_h^{1,p}(V)}\le {h_0}^{\frac{1}{p}}\left(\frac{\lambda_1\|e_1\|_{L^{\frac{p}{p-1}}(V)}+\|g_1\|_{L^{\frac{p}{p-1}}(V)}}
{h_0-\|f_1\|_{\infty}}\right)^{\frac{1}{p-1}},
$$
then by (\ref{bb1}), we have
$$
\|u_{\lambda \star}\|_{\infty}\le {\mu_0}^{-\frac{1}{p}}\left(\frac{\lambda_1\|e_1\|_{L^{\frac{p}{p-1}}(V)}+\|g_1\|_{L^{\frac{p}{p-1}}(V)}}
{h_0-\|f_1\|_{\infty}}\right)^{\frac{1}{p-1}}.
$$
Similarly, when $(u_{\lambda \star},v_{\lambda \star})=(0,v_{\lambda \star})$, we have
$$
\|v_{\lambda \star}\|_{\infty}\le {\mu_0}^{-\frac{1}{q}}\left(\frac{\lambda_2\|e_2\|_{L^{\frac{q}{q-1}}(V)}+\|g_2\|_{L^{\frac{q}{q-1}}(V)}}
{h_0-\|f_2\|_{\infty}}\right)^{\frac{1}{q-1}}.
$$
\qed

\vskip2mm
\noindent
{\bf Proof of Theorem 1.2.} The proof is similar to Theorem 1.1, in which we only need to slightly modify the proof of Lemma 3.1
with replacing $(F_1)$ by $(F_1')$. We omit the details. \qed

\vskip2mm
\section{Proofs for the super-$(p,q)$-linear case}
\setcounter{equation}{0}
\noindent
{\bf Lemma 4.1.} {\it Assume that $(H_1)$ and $(C_1)$ hold.
Then for each $\lambda\in(0,{\lambda_0})$, there exists a positive constant $\rho_\lambda$ such that $\varphi(u,v)>0$ whenever $||(u,v)||_W=\rho_\lambda$.}\\
{\bf Proof.}  Note that $F(x,0,0)=0$. By $(C_1)$, for all $(s,t)\in \R^2$ with $|(s,t)|< l_0$, we have
\begin{eqnarray}
\label{EQ8}
\notag |F(x,s,t)|&=&|F(x,s,t)-F(x,0,0)|\\
          \notag &\le& |F(x,s,t)-F(x,0,t)|+|F(x,0,t)-F(x,0,0)|\\
          \notag &\le& \int_0^{|s|}|F_s(x,s,t)|ds+\int_0^{|t|}|F_t(x,0,t)|dt\\
          \notag &\le& \int_0^{|s|}\frac{h_0}{q+1}\left(|s|^{p -1}+|t|^{\frac{pq-q}{p}}\right)ds+\int_0^{|t|}\frac{h_0}{q+1}|t|^{q -1} dt\\
          \notag &\le& \frac{h_0}{p(q+1)}|s|^{p}+\frac{h_0}{q+1}|t|^{\frac{pq-q}{p}}|s|+\frac{h_0}{q(q+1)}|t|^{q}\\
                 &\le& \frac{2h_0}{p(q+1)}|s|^{p }+\frac{(pq-q+p)h_0}{pq(q+1)}|t|^{q}.
\end{eqnarray}
 It is easy to obtain that for each $\lambda$ satisfying (\ref{pp1}), there exists a $\varepsilon_\lambda>0$ such that
$$
0<\lambda<\lambda_\varepsilon:=\frac{\min\{1,q-1\}}{2^{\max\{p,q\}-1}(pq+p)\max\left\{{h_0}^{-\frac{1}{p}}\|e_1\|_{L^{\frac{p}{p-1}}(V)}
,{h_0}^{-\frac{1}{q}}\|e_2\|_{L^{\frac{q}{q-1}}(V)}\right\}}\cdot\left(\Lambda_0-\varepsilon_\lambda\right)^{\max\{p,q\}-1}.
$$
For any $(u,v)\in W$ with $\|(u,v)\|_W=\Lambda_0-\varepsilon_\lambda$,
by  (\ref{bb1}), we have $\|u\|_\infty< \frac{l_0}{2}$ and $\|v\|_\infty< \frac{l_0}{2}$, and so $|(u(x),v(x))|\le \|u\|_\infty+\|v\|_\infty<l_0$
for all $x\in V$. Then
\begin{eqnarray}
 \label{pp7}
 &       & \varphi_\lambda(u,v)\nonumber\\
 & \ge  & \frac{1}{p}\|u\|_{W_h^{1,p}(V)}^p+\frac{1}{q}\|v\|_{W_h^{1,q}(V)}^q- \frac{2h_0}{p(q+1)}\int_V|u|^pd\mu-\frac{(pq-q+p)h_0}{pq(q+1)}\int_V|v|^q d\mu\nonumber\\
 &       &  -\lambda\int_V(e_1u+e_2v)d\mu\nonumber\\
 & \ge  &  \left(\frac{1}{p}-\frac{2}{p(q+1)}\right)\|u\|_{W_h^{1,p}(V)}^p+\left(\frac{1}{q}-\frac{pq-q+p}{pq(q+1)}\right)\|v\|_{W_h^{1,q}(V)}^q
            -\lambda\int_V(e_1u+e_2v)d\mu\nonumber\\
 & \ge  & \frac{\min\{1,q-1\}}{pq+p}\left(\|u\|_{W_h^{1,p}(V)}^p+\|v\|_{W_h^{1,q}(V)}^q\right)
            -\lambda\max\left\{{h_0}^{-\frac{1}{p}}\|e_1\|_{L^{\frac{p}{p-1}}(V)},{h_0}^{-\frac{1}{q}}\|e_2\|_{L^{\frac{q}{q-1}}(V)}\right\}\|(u,v)\|_W\nonumber\\
 & \ge  &  \frac{\min\{1,q-1\}}{2^{\max\{p,q\}-1}(pq+p)}\|(u,v)\|_W^{\max\{p,q\}}
           -\lambda\max\left\{{h_0}^{-\frac{1}{p}}\|e_1\|_{L^{\frac{p}{p-1}}(V)},{h_0}^{-\frac{1}{q}}\|e_2\|_{L^{\frac{q}{q-1}}(V)}\right\}\|(u,v)\|_W
\end{eqnarray}
for any $(u,v)\in W$ with $\|(u,v)\|_W=\Lambda_0-\varepsilon_\lambda$. Let $\rho_\lambda=\Lambda_0-\varepsilon_\lambda$.
Hence, for each $\lambda\in(0,{\lambda_0})$, there exists a  $\rho_\lambda$ such that $\varphi(u,v)\ge \alpha_{\lambda}>0$ whenever $||(u,v)||_W=\rho_\lambda$, where
\begin{eqnarray}
\label{pp8}\alpha_{\lambda}=  \frac{\min\{1,q-1\}}{2^{\max\{p,q\}-1}(pq+p)}\rho_\lambda^{\max\{p,q\}}
           -\lambda\max\left\{{h_0}^{-\frac{1}{p}}\|e_1\|_{L^{\frac{p}{p-1}}(V)},{h_0}^{-\frac{1}{q}}\|e_2\|_{L^{\frac{q}{q-1}}(V)}\right\}\rho_\lambda.
\end{eqnarray}
\qed

\vskip2mm
\noindent
{\bf Lemma 4.2.} {\it Assume that $(C_2)$ holds. Then for each $\lambda\in(0,{\lambda_0})$,
there exists a $(u_{**\lambda},v_{**\lambda})\in W$ with $\|(u_{**\lambda},v_{**\lambda})\|_W>\rho_\lambda$ such that $\varphi(u_{**\lambda},v_{**\lambda})<0$.}\\
{\bf Proof.} Let
$$
u^*(x)=v^*(x)=\begin{cases}
1,& x= x_3,\\
0,& x\not=x_3,
\end{cases}
$$
where $x_3 \in V$ with $\mu(x_3)>0$ and $e_1(x_3)+e_2(x_3)>0$. Then
\begin{eqnarray}\label{M1}
          \int_V |\nabla u^*|^pd\mu
&   =   & \sum_{x\in V}|\nabla u^*|^p(x)\mu(x)\nonumber\\
&   =   & \sum_{x\in V}\left(\frac{1}{2\mu(x)}\sum\limits_{y\thicksim x}w_{xy}(u^*(y)-u^*(x))^2\right)^{\frac{p}{2}}\mu(x)\nonumber\\
&   =   & \left(\frac{1}{2\mu(x_3)}\sum\limits_{y\thicksim x_3}w_{x_3y}\right)^{\frac{p}{2}}\mu(x_3)+\sum_{x\thicksim x_3}\left(\frac{1}{2\mu(x)}\sum\limits_{x_3\thicksim x}w_{xx_3}\right)^{\frac{p}{2}}\mu(x)\nonumber\\
&   =   & \left(\frac{deg(x_3)}{2\mu(x_3)}\right)^{\frac{p}{2}}\mu(x_3)+\sum_{x\thicksim x_3}\left(\frac{deg(x_3)}{2\mu(x)}\right)^{\frac{p}{2}}\mu(x)\nonumber\\
&   =   & \left(\frac{deg(x_3)}{2}\right)^{\frac{p}{2}}\left(\sum_{x\thicksim x_3}\left(\frac{1}{\mu(x)}\right)^{\frac{p}{2}-1}+\frac{1}{\mu(x_3)^{\frac{p}{2}-1}}\right)\nonumber\\
&  :=   & D_1.
\end{eqnarray}
Similarly, we have
$$
 \int_V |\nabla v^*|^qd\mu=\left(\frac{deg(x_3)}{2}\right)^{\frac{q}{2}}\left(\sum_{x\thicksim x_3}\left(\frac{1}{\mu(x)}\right)^{\frac{q}{2}-1}+\frac{1}{\mu(x_3)^{\frac{q}{2}-1}}\right):=D_2.
$$
Thus, by $(C_2)$, for all $s\in \R$ with $s>l_1$,  we have
\begin{eqnarray}
\label{EQ6}
           \notag \varphi_\lambda(su^*,sv^*)
&   =   & \frac{s^p}{p}\|u^*\|_{W_h^{1,p}(V)}^p+\frac{s^q}{q}\|v^*\|_{W_h^{1,q}(V)}^q-\int_VF(x,su^*(x),sv^*(x))d\mu\nonumber\\
&       &  -\lambda \int_V(se_1u^*+se_2v^*)d\mu\nonumber\\
&   =   & \frac{s^p}{p}\left(D_1+\mu(x_3)h_1(x_3)\right)+\frac{s^q}{q}\left(D_2+\mu(x_3)h_2(x_3)\right)\nonumber\\
&       &  -\mu(x_3)F(x_3,s,s)-\lambda s\mu(x_3)(e_1(x_3)+e_2(x_3))\nonumber\\
& \le   & \frac{s^p}{p}\left(D_1+\mu(x_3)h_1(x_3)\right) +\frac{s^q}{q}\left(D_2+\mu(x_3)h_2(x_3)\right)\nonumber\\
&       &  -M\mu(x_3)(s^p+s^q)-\lambda s\mu(x_3)(e_1(x_3)+e_2(x_3))\nonumber\\
&   =   &  s^p\left(\frac{D_1+\mu(x_3)h_1(x_3)}{p}-M\mu(x_3)\right) +s^q\left(\frac{D_2+\mu(x_3)h_2(x_3)}{q}-M\mu(x_3)\right)\nonumber\\
&       &  -\lambda s\mu(x_3)(e_1(x_3)+e_2(x_3)),
\end{eqnarray}
which implies $\varphi(su^*,sv^*)\rightarrow -\infty$ as $s\rightarrow +\infty$. Hence, for each $\lambda\in(0,{\lambda_0})$,
 there exists $s_\lambda$ large enough such that  $\|(s_\lambda u^*,s_\lambda v^*)\|_W>\rho_\lambda$  and  $\varphi(s_\lambda u^*,s_\lambda v^*)<0$.
 Let $u_{**\lambda}=s_\lambda u^*$ and $v_{**\lambda}=s_\lambda v^*$. Then the proof is completed. \qed

\vskip2mm
\noindent
{\bf Lemma 4.3.} {\it Assume that $(F_0), (C_3)$, $(H_1)$ and $(H_2')$ hold.
  Then for each $\lambda\in(0,{\lambda_0})$, $\varphi_\lambda$ satisfies the (PS)-condition.
 } \\
{\bf Proof.} Let $\{(u_k,v_k)\} \subset W$ be a Palais-Smale sequence of $\varphi_{\lambda}$. Then there exists a positive constant $c$ such that
$$
|\varphi_\lambda(u_k,v_k)|\le c \mbox{ for all }k\in \mathbb N \ \ \mbox{and} \;\;\;\varphi_\lambda'(u_k,v_k)\rightarrow 0 \ \ \mbox{as } k\to\infty.
$$
Then, by $(C_3)$, we have
\begin{eqnarray}
\label{pp2}&       &   c+\|u_k\|_{W_h^{1,p}(V)}+\|v_k\|_{W_h^{1,q}(V)}\nonumber\\
&   =   &             c+{\|(u_k,v_k)\|}_W\nonumber\\
& \ge  & \varphi_\lambda(u_k,v_k)-\frac{1}{\nu}\langle\varphi_\lambda'(u_k,v_k),(u_k,v_k)\rangle\nonumber\\
&  =    & \left(\frac{1}{p}-\frac{1}{\nu}\right)\|u_k\|_{W_h^{1,p}(V)}^p+\left(\frac{1}{q}-\frac{1}{\nu}\right)\|{v_k}\|_{W_h^{1,q}(V)}^q\nonumber\\
&       & -\frac{1}{\nu}\int_V\left[\nu F(x,u_k,v_k)-F_u(x,u_k,v_k)u_k-F_v(x,u_k,v_k)v_k\right]d\mu\nonumber\\
&       & -\frac{\nu-1}{\nu}\lambda \int_V(e_1u_k+e_2v_k)d\mu\nonumber\\
& \ge  & \left(\frac{1}{p}-\frac{1}{\nu}\right)\|u_k\|_{W_h^{1,p}(V)}^p+\left(\frac{1}{q}-\frac{1}{\nu}\right)\|{v_k}\|_{W_h^{1,q}(V)}^q\nonumber\\
&       & -\frac{A}{\nu}\int_V( |u|^p+|v|^q) d\mu-\frac{\nu-1}{\nu}\lambda \int_V(e_1u_k+e_2v_k)d\mu\nonumber\\
& \ge  & \left(\frac{1}{p}-\frac{1}{\nu}-\frac{A}{\nu h_0}\right)\|u_k\|_{W_h^{1,p}(V)}^p+\left(\frac{1}{q}-\frac{1}{\nu}-\frac{A}{\nu h_0}\right)\|{v_k}\|_{W_h^{1,q}(V)}^q\nonumber\\
&       & -\frac{(\nu-1)\lambda}{\nu} \left(h_0^{-\frac{1}{p}}\|e_1\|_{L^{\frac{p}{p-1}}(V)}\|u_k\|_{W_h^{1,p}(V)}
          +h_0^{-\frac{1}{q}}\|e_2\|_{L^{\frac{q}{q-1}}(V)}\|v_k\|_{W_h^{1,q}(V)}\right).
\end{eqnarray}
We claim that ${\|(u_k,v_k)\|}_W$ is bounded. In fact, if
\begin{eqnarray}\label{qq1}
\|u_k\|_{W_h^{1,p}(V)}\to\infty  \mbox{ and } \|v_k\|_{W_h^{1,q}(V)}\to \infty, \mbox{ as } k\to \infty,
\end{eqnarray}
then (\ref{pp2}) implies that
\begin{eqnarray}
\label{pp3}&       &  c+{\|(u_k,v_k)\|}_W +\frac{(\nu-1)\lambda}{\nu}\max\left\{h_0^{-\frac{1}{p}}\|e_1\|_{L^{\frac{p}{p-1}}(V)},h_0^{-\frac{1}{q}}\|e_2\|_{L^{\frac{q}{q-1}}(V)}\right\}
          {\|(u_k,v_k)\|}_W\nonumber\\
& \ge  & \min\left\{\left(\frac{1}{p}-\frac{1}{\nu}-\frac{A}{\nu h_0}\right),\left(\frac{1}{q}-\frac{1}{\nu}-\frac{A}{\nu h_0}\right)\right\}\left(\|u_k\|_{W_h^{1,p}(V)}^p+\|{v_k}\|_{W_h^{1,q}(V)}^q\right)\nonumber\\
& \ge  & \min\left\{\left(\frac{1}{p}-\frac{1}{\nu}-\frac{A}{\nu h_0}\right),\left(\frac{1}{q}-\frac{1}{\nu}-\frac{A}{\nu h_0}\right)\right\}\frac{1}{2^{\min\{p,q\}-1}}}\|(u_k,v_k)\|_W^{\min\{p,q\}
\end{eqnarray}
for all large $k$, which contradicts with (\ref{qq1}).
If
\begin{eqnarray}\label{qq2}
\|u_k\|_{W_h^{1,p}(V)}\to\infty \mbox{ as } k\to \infty
\end{eqnarray}
and $\|v_k\|_{W_h^{1,q}(V)}$ is bounded for all $k\in \mathbb N$, then by (\ref{pp2}), there exists
two positive constants $c_0$ and $c_1$ such that
\begin{eqnarray*}
             c_0+c_1\|u_k\|_{W_h^{1,p}(V)}
 \ge   \left(\frac{1}{p}-\frac{1}{\nu}-\frac{A}{\nu h_0}\right)\|u_k\|_{W_h^{1,p}(V)}^p
\end{eqnarray*}
which contradicts with (\ref{qq2}). Similarly,
if $\|v_k\|_{W_h^{1,q}(V)}\to\infty$ as $k\to \infty$ and $\|u_k\|_{W_h^{1,p}(V)}$ is bounded for all $k\in \mathbb N$,
we can also obtain the same contradiction. Hence, the above arguments imply that both $\|u_k\|_{W_h^{1,p}(V)}$ and $\|v_k\|_{W_h^{1,q}(V)}$
are bounded. So there exists a positive constant $c_2$ such that $\|u_k\|_{W_h^{1,p}(V)}\le c_2$ and $\|v_k\|_{W_h^{1,q}(V)}\le c_2$.
Then we can find a
subsequence, still denoted by $\{u_k\}$, such that $u_k\rightharpoonup u^{\star}_{\lambda}$ for some $u^{\star}_{\lambda}\in W_h^{1,p}(V)$ as $k\to \infty$, and a
subsequence of $\{v_k\}$,  which has the same subscript as the subsequence of $\{u_k\}$, still denoted by $\{v_k\}$,
such that $v_k\rightharpoonup v^{\star}_{\lambda}$ for some $v^{\star}_{\lambda}\in W_h^{1,q}(V)$ as $k\to \infty$.
By Lemma 2.2, we know that
\begin{eqnarray}
\label{pp4}&   &  u_k\to  u^{\star}_{\lambda} \mbox{ and } v_k\to v^{\star}_{\lambda} \ \mbox{in } L^{\infty}(V), \ \ \mbox{as }k\to \infty.
\end{eqnarray}
   Then by (\ref{eq3}), we have
\begin{eqnarray*}
\label{pp6}
\notag &    & \langle\varphi_\lambda'(u_k,v_k)-\varphi_\lambda'(u^{\star}_{\lambda},v^{\star}_{\lambda}),(u_k-u^{\star}_{\lambda},0)\rangle\\
\notag & =  & \|u_k\|_{W_h^{1,p}(V)}^p+\|u^{\star}_{\lambda}\|_{W_h^{1,p}(V)}^p-\int_V \left[|\nabla u_k|^{p-2}\Gamma( u_k,u^{\star}_{\lambda})+h_1(x)|u_k|^{p-2}u_ku^{\star}_{\lambda}\right]d\mu\\
       &    & -\int_V \left[|\nabla u^{\star}_{\lambda}|^{p-2}\Gamma(u^{\star}_{\lambda},u_k)+h_1(x)|u^{\star}_{\lambda}|^{p-2}u^{\star}_{\lambda}u_k \right]d \mu+\int_V[F_u(x,u^{\star}_{\lambda},v^{\star}_{\lambda})-F_u(x,u_k,v_k)](u_k-u^{\star}_{\lambda})d\mu.
\end{eqnarray*}
Let $ A_1=c_2 \frac{1}{h_0^{1/p}\mu_0^{1/p}}+c_2 \frac{1}{h_0^{1/q}\mu_0^{1/q}}$
and $ A_2=\|u^{\star}_{\lambda}\|_\infty+\|v^{\star}_{\lambda}\|_\infty$.  By $(F_0)$ and (\ref{pp4}), we have
\begin{eqnarray*}
\notag &       & \int_V[F_u(x,u^{\star}_{\lambda},v^{\star}_{\lambda})-F_u(x,u_k,v_k)](u_k-u^{\star}_{\lambda})d\mu\\
\notag & \le  & \int_V|F_u(x,u^{\star}_{\lambda},v^{\star}_{\lambda})-F_u(x,u_k,v_k)||u_k-u^{\star}_{\lambda}|d\mu\\
\notag & \le  & \int_V[|F_u(x,u_k,v_k)|+|F_u(x,u^{\star}_{\lambda},v^{\star}_{\lambda})|]|u_k-u^{\star}_{\lambda}|d\mu\\
\notag & \le  & \left[\max_{|(s,t)|\le A_1}a(|(s,t)|) \int_V b(x)d\mu+ \max_{|(s,t)|\le A_2}a(|(s,t)|) \int_V b(x)d\mu\right] \|u_k-u^{\star}_{\lambda}\|_\infty \\
       &\rightarrow& 0.
\end{eqnarray*}
The rest of arguments are the same as Lemma 3.2.
\qed

\vskip2mm
 \noindent
 {\bf Lemma 4.4.} {\it Assume that $(C_1)$ and $(C_4)$ holds. Then for each $\lambda\in (0,\lambda_0)$, $-\infty<\inf\{\varphi(u,v):(u,v)\in\bar B_{\rho_\lambda} \}<0$,
 where $\rho_\lambda$ is given in Lemma 4.1 and $\bar B_{\rho_\lambda}=\{(u,v)\in W\big|\|(u,v)\|_W\le \rho_\lambda\}$}. \\
{\bf Proof.} Let
$$
u^{**}(x)=v^{**}(x)=\begin{cases}
1,& x= x_4,\\
0,& x\not=x_4,
\end{cases}
$$
where $x_4 \in V$ with $\mu(x_4)>0$ and $e_1(x_4)+e_2(x_4)>0$. Hence, by (\ref{M1}), we obtain that
\begin{eqnarray*}
D_3:=\int_V |\nabla u^{**}|^pd\mu=\left(\frac{deg(x_4)}{2}\right)^{\frac{p}{2}}\left(\sum_{x\thicksim x_4}\left(\frac{1}{\mu(x)}\right)^{\frac{p}{2}-1}+\frac{1}{\mu(x_4)^{\frac{p}{2}-1}}\right),\\
D_4:=\int_V |\nabla v^{**}|^qd\mu=\left(\frac{deg(x_4)}{2}\right)^{\frac{q}{2}}\left(\sum_{x\thicksim x_4}\left(\frac{1}{\mu(x)}\right)^{\frac{q}{2}-1}+\frac{1}{\mu(x_4)^{\frac{q}{2}-1}}\right).
\end{eqnarray*}
Then for each $\lambda\in (0,\lambda_0)$,  by $(C_4)$, for all $t\in \R$ with $0<t<l_2$,  we have
\begin{eqnarray}
\label{EQ6}
           \notag \varphi_\lambda(tu^{**},tv^{**})
&   =   & \frac{t^p}{p}\|u^{**}\|_{W_h^{1,p}(V)}^p+\frac{t^q}{q}\|v^{**}\|_{W_h^{1,q}(V)}^q-\int_VF(x,tu^{**}(x),tv^{**}(x))d\mu\nonumber\\
&       &  -\lambda \int_V(te_1u^{**}+te_2v^{**})d\mu\nonumber\\
&   =   & \frac{t^p}{p}\left(D_3+\mu(x_4)h_1(x_4)\right)+\frac{t^q}{q}\left(D_4+\mu(x_4)h_2(x_4)\right)\nonumber\\
&       &  -\mu(x_4)F(x_4,t,t)-\lambda t\mu(x_4)\left(e_1(x_4)+e_2(x_4)\right)\nonumber\\
& \le  & \frac{t^p}{p}\left(D_3+\mu(x_4)h_1(x_4)\right)+\frac{t^q}{q}\left(D_4+\mu(x_4)h_2(x_4)\right)\nonumber\\
&       &  +K_3(x_4)\mu(x_4)|t|^{\beta_3}-\lambda t\mu(x_4)\left(e_1(x_4)+e_2(x_4)\right).
\end{eqnarray}
 Note that $p>1$, $q>1$, $\beta_3>1$ and $K_3(x_4)>0$. By $(\ref {EQ6})$, there exists a  sufficiently small $t_{1,\lambda}$ satisfying
 $$
 0<t_{1,\lambda}<\min\left\{\frac{\rho_\lambda}{2\|u^{**}\|_{W_h^{1,p}(V)}},\frac{\rho_\lambda}{2\|v^{**}\|_{W_h^{1,q}(V)}}\right\}
 $$
 such that $\varphi(t_{1,\lambda}u^{**},t_{1,\lambda}v^{**}) < 0$.
 Obviously,  $\|(t_{1,\lambda}u^{**},t_{1,\lambda}v^{**})\|_{W}<\rho_\lambda$.
 Hence, $\inf\{\varphi(u,v):(u,v)\in\bar B_{\rho_\lambda} \} \le \varphi(t_{1,\lambda}u^{**},t_{1,\lambda}v^{**})<0$. Moreover, it is easy to see that
 (\ref{pp7}) still holds for all $(u,v)\in\bar B_{\rho_\lambda}$. Then
 \begin{eqnarray*}
        \varphi_\lambda(u,v)
 \ge    -\lambda\max\left\{{h_0}^{-\frac{1}{p}}\|e_1\|_{L^{\frac{p}{p-1}}(V)},{h_0}^{-\frac{1}{q}}\|e_2\|_{L^{\frac{q}{q-1}}(V)}\right\}\rho_\lambda,
\end{eqnarray*}
which shows that $\varphi_\lambda$ is bounded from below in $\bar B_{\rho_\lambda}$ for each $\lambda\in (0,\lambda_0)$. So $\inf\{\varphi(u,v):(u,v)\in\bar B_{\rho_\lambda}\}>-\infty$.
 \qed

 \vskip2mm
\noindent
{\bf Proof of Theorem 1.3.} By Lemma 4.1-Lemma 4.3 and Lemma 2.5, we obtain that for each $\lambda\in(0,{\lambda_0})$,
  $\varphi_\lambda$ has a critical value $c_*\ge \alpha_{\lambda}>0$ with
 $$
 c_*:=\inf_{\gamma\in\Gamma}\max_{t\in[0,1]}\varphi_\lambda(\gamma(t)),
 $$
where
 $$
 \Gamma:=\{\gamma\in C([0,1],X]):\gamma(0)=(0,0),\gamma(1)=(u_{*,\lambda},v_{*,\lambda})\}
 $$
 and $\alpha_{\lambda}$ is defined by (\ref{pp8}). Hence, by Proposition 3.1,  system (\ref{eq1}) has one  solution $(u^{\lambda \star},v^{\lambda \star})$ of  positive energy.
 Obviously, $(u^{\lambda \star},v^{\lambda \star})\not=(0,0)$. Otherwise, by the fact that $F(x,0,0)=0$ for all $x\in V$,
 we have $\varphi(u^{\lambda \star},v^{\lambda \star})=0$, which contradicts with $c_*>0$.
 \par
 Next, we prove that system (\ref{eq1}) has one  solution  of  negative energy if $(C_4)$ also holds.
 The proof is motivated by Theorem 3.3 in \cite{Cheng}. In fact,
 by Lemma 4.1 and Lemma 4.4, we know that
$$
-\infty<\inf_{\bar B_{\rho_\lambda}}\varphi_{\lambda}<0<\inf_{\partial B_{\rho_\lambda}}\varphi_{\lambda}
$$
for each $\lambda \in (0,\lambda_0)$. Set
$$
\frac{1}{n}\in \left(0,\inf_{\partial B_{\rho_\lambda}}\varphi_{\lambda}-\inf_{\bar B_{\rho_\lambda}} \varphi_{\lambda}\right),\;n\in \Z^+.
$$
Then there exists a $(u_n,v_n)\in \bar B_{\rho_\lambda}$ such that
\begin{eqnarray}
\label{EQ13}
\varphi_{\lambda}(u_n,v_n)\le\inf_{\bar B_{\rho_\lambda}}\varphi_{\lambda}+\frac {1}{n}.
\end{eqnarray}
As $\varphi_{\lambda}(u,v) \in C^1(W,\R)$, we know $\varphi_{\lambda}(u,v)$ is lower semicontinuous. Thus, by Lemma 2.3 we have
$$
\varphi_{\lambda}(u_n,v_n)\le \varphi_{\lambda}(u,v)+\frac{1}{n}\|(u,v)-(u_n,v_n)\|_W,\;\forall (u,v)\in \bar B_{\rho_\lambda}.
$$
Note that
$$
\varphi_{\lambda}(u_n,v_n)\le\inf_{\bar B_{\rho_\lambda}}\varphi_{\lambda}+\frac {1}{n}<\inf_{\partial B_{\rho_\lambda}}\varphi_{\lambda}.
$$
Thus, $(u_n,v_n)\in B_{\rho_\lambda}$. Defining $ M_n:W \rightarrow R$ by
$$
M_n(u,v)=\varphi_{\lambda}(u,v) + \frac {1}{n}\|(u,v)-(u_n,v_n)\|_W,
$$
we have $(u_n,v_n)\in B_{\rho_\lambda}$ minimizes $M_n$ on $ \bar B_{\rho_\lambda}$. Therefore, for all $(u,v)\in W$ with $\|(u,v)\|_W = 1$,
taking $t >0$ small enough such that $(u_n+tu,v_n+tv)\in \bar B_{\rho_\lambda}$, then
$$
\frac{M_n(u_n+tu,v_n+tv)-M_n(u_n,v_n)}{t}\ge 0,
$$
which implies that
$$
\langle\varphi_{\lambda}'(u_n,v_n),(u,v)\rangle\ge-\frac{1}{n}.
$$
Similarly, when $t<0$ and $|t|$ small enough, we have
$$
\langle\varphi_{\lambda}'(u_n,v_n),(u,v)\rangle\le\frac{1}{n}.
$$
Hence,
\begin{eqnarray}
\label{EQ14}
\|\varphi_{\lambda}'(u_n,v_n)\|\le\frac{1}{n}.
\end{eqnarray}
Passing to the limit in (\ref{EQ13}) and (\ref{EQ14}), we conclude that
$$
\varphi_{\lambda}(u_n,v_n)\rightarrow \inf_{\bar B_{\rho_\lambda}} \varphi_{\lambda}\;\;\;\mbox{and} \;\;\;\|\varphi_{\lambda}'(u_n,v_n)\|\rightarrow 0\;\;\;as \;n\rightarrow\infty.
$$
Hence, $\{(u_n,v_n)\}\subset \bar B_{\rho_\lambda} $ is a Palais-Smale sequence of $\varphi_{\lambda}$.
By Lemma 4.3, $\{(u_n,v_n)\}$ has a strongly convergent subsequence $\{(u_{nk},v_{nk})\}\subset \bar B_{\rho_\lambda}$,
 and $(u_{nk},v_{nk})\rightarrow(u^{\star\star},v^{\star\star})\in \bar B_{\rho_\lambda}$ as $n_k \rightarrow\infty$. Consequently,
$$
\varphi_{\lambda}(u^{\star\star},v^{\star\star})=\inf_{\bar B_{\rho_\lambda}}\varphi_{\lambda}<0\;\;\;\mbox{and} \;\;\;\varphi_{\lambda}'(u^{\star\star},v^{\star\star})=0,
$$
which implies that  system (\ref{eq1})  has  a solution $(u^{\star\star},v^{\star\star})\not=(0,0)$ of negative energy.
\qed

\vskip2mm
\section{Examples}
\setcounter{equation}{0}
\vskip2mm
\noindent
{\bf Example 5.1.} Let $p=2$ and $q=3$. Consider the following system:
\begin{eqnarray}
\label{eaxmple 1}
\begin{cases}
-\Delta u+h_1(x)u=F_u(x,u,v)+\lambda_1 e_1(x),x\in V,\\
-\Delta_3 v+h_2(x)v=F_v(x,u,v)+\lambda_2 e_2(x),x \in V,
\end{cases}
\end{eqnarray}
where $G=(V,E)$ is locally finite graph, the measure $\mu(x)\ge \mu_0=1$ for all $x \in V$, $h_i:V \to \R^+,\;i=1,2$, $h_1(x)=3+dist(x,x_1),\;h_2(x)=3+dist(x,x_2)$, where $x_1$ and $x_2$
are two fixed points in $V$ and $\mu(x_1)=\mu(x_2)=1$,
\begin{eqnarray*}
F(x,s,t) & = &
\begin{cases}
\frac{3}{5}\left(s^{\frac{5}{3}}+t^{\frac{5}{3}}\right)&,x=x_1,x_2,\\
\ \ 0&,x\neq x_1,x_2,
\end{cases}\\
e_1(x)= e_2(x)& = &
\begin{cases}
1&,x=x_1,x_2\\
0&,x\neq x_1,x_2,
\end{cases}
\end{eqnarray*}
and $\lambda_1,\lambda_2>0$.
Next, We verify that $h_1,h_2$ and $F$ satisfy the conditions in Theorem 1.1:
\begin{itemize}
  \item  Obviously, when $dist(x,x_i)\rightarrow +\infty,h_i(x)\rightarrow +\infty$ and $h_i\ge h_0=3,i=1,2$. \\
         Hence, $h_i$ satisfies $(H_1),(H_2)$, $i=1,2$.
  \item  Let
         $$
         f_1(x)\equiv1,\;
         g_1(x)=
         \begin{cases}
         1&,x=x_1,x_2,\\
         0&,x\neq x_1,x_2.
         \end{cases}
         $$
         Then
         $$
          \|f_1\|_{\infty}=1<\min\left\{\frac{h_0}{2},\frac{ph_0}{q(p-1)}\right\}=\frac{3}{2},\ \
          \|g_1\|_{L^2(V)}=\sqrt{2}.
         $$
         Moreover,
         $$
         |F_s(x,s,t)|
         =   |s|^{\frac{2}{3}}
         \le |s|+1
         \le f_1(x)(|s|+|t|^{\frac{3}{2}})+g_1(x).
         $$
         Similarly, let
         $$
         f_2(x)=f_1(x),\;g_2(x)=g_1(x).
         $$
         We also have
         $$
         |F_t(x,s,t)|
         =   |t|^{\frac{2}{3}}
         \le |t|^2+1
         \le f_2(x)(|s|^2+|t|^2)+g_2(x).
         $$
         Then $F(x,s,t)$ satisfies $(F_1)$. Hence, $F(x,s,t)$ also satisfies $(F_0)$ .
  \item  Let
         $$
         \beta_1=\frac{5}{3},K_1(x)\equiv\frac{3}{5}.
         $$
         Then
         $$
         F(x,s,0)\ge-\frac{3}{5}|s|^{\frac{5}{3}}.
         $$
         Hence, $F(x,s,t)$ satisfies (i) of $(F_2)$.
\end{itemize}
Hence, by Theorem 1.1, for each pair $(\lambda_1,\lambda_2) \in (0,+\infty)\times(0,+\infty)$, system (\ref{eaxmple 1}) has one nontrivial solution $(u_{\lambda \star},v_{\lambda \star})$. Furthermore, if $(u_{\lambda \star},v_{\lambda \star})=(u_{\lambda \star},0)$, then
 $\|u_{\lambda \star}\|_{\infty}\le \frac{\sqrt{2}}{2}(\lambda_1+1)$. If $(u_{\lambda \star},v_{\lambda \star})=(0,v_{\lambda \star})$, then
$\|v_{\lambda \star}\|_{\infty}\le 2^{\frac{1}{3}}\left(\frac{\lambda_2+1}{2}\right)^{\frac{1}{2}}$.

\vskip2mm
\noindent
{\bf Example 5.2.} Let $p=2$ and $q=3$. Consider the following system:
\begin{eqnarray}
 \label{eaxmple 2}
 \begin{cases}
  -\Delta u+h_1(x)u=F_u(x,u,v)+\lambda e_1(x),x\in V,\\
  -\Delta_3 v+h_2(x)v=F_v(x,u,v)+\lambda e_2(x),x \in V,
 \end{cases}
\end{eqnarray}
where $G=(V,E)$ is locally finite graph, the measure $\mu(x)\ge \mu_0> 0$ for all $x \in V$, $e_1,e_2\in L^2(V),e_1(x),\;e_2(x)\not\equiv 0$ and $\lambda>0$.
$h_1(x)=h_2(x)=c_1dist(x,x_1)-\frac{1}{dist(x,x_2)+1}+2$, where $c_1$ is positive constant and $x_1$ and $x_2$
are two fixed points in $V$ with $e_1(x_1)+e_2(x)>0$, $F(x,s,t)=M\ln(1+s^4+t^4)(s^4+t^4)$,
$$M=\max\left\{\frac{D_1+\mu(x_1)h_1(x_1)}{2\mu(x_1)},
\frac{D_2+\mu(x_1)h_2(x_1)}{3\mu(x_1)}\right\}+1>1,$$
 $D_1=\frac{deg(x_1)}{2}\left(\sharp A+1\right)$, and\\
  $$D_2=\left(\frac{deg(x_1)}{2}\right)^{\frac{3}{2}}\left(\sum_{x\thicksim x_1 }\frac{1}{\sqrt{\mu(x)}}+\frac{1}{\sqrt{\mu(x_1)}}\right),$$
  where $\sharp A$ is the number of elements in the set $A=\{x\in V \big| x\thicksim x_1\}$.

\par
We verify that $h_1,h_2$ and $F$ satisfy the conditions in Theorem 1.3:
\begin{itemize}
 \item  Obviously, $h_1,h_2$ satisfy $(H_1)$ and $h_0=1$.
 \item  For any given constant $B$, when $h_1=h_2=c_1dist(x,x_1)-\frac{1}{dist(x,x_2)+1}+2<B$, we have
        $$
        c_1dist(x,x_1)<B-2+\frac{1}{dist(x,x_2)+1}<B-1.
        $$
        Moreover, since $V$ is a locally finite graph, the set $A_i=\{x\in V\big|h_i \le B\}\subseteq\{x\in V\big||x-x_1|<B-1\}$ is finite. So, $\sum_{x \in A_i}\mu(x)$ is finite,
         $(H_2')$ holds.
 \item  By $F(x,s,t)$, when $|s|$ and $|t| <1$, we have
       \begin{eqnarray*}
        |F_s(x,s,t)|
        &  =  & M\left|\frac{4s^3(s^4+t^4)}{1+s^4+t^4}+4s^3\ln(1+s^4+t^4)\right|\\
        & \le & 4M\left(\frac{|s|^3(s^4+t^4)}{1+s^4+t^4}+|s|^3(1+s^4+t^4)\right)\\
        & \le & 4M(4|s|^3+t^4).
       \end{eqnarray*}
       Moreover, when $|s|<\frac{1}{16\sqrt{M}}$ and $|t|<\frac{1}{(16M)^{\frac{2}{5}}}$, we have
        $$
        16M|s|^3\le \frac{1}{4}|s|,4Mt^4\le \frac{1}{4}|t|^{\frac{3}{2}}.
        $$
        So, when $|(s,t)|\le \frac{1}{16\sqrt{M}}$, $F_s(x,s,t)|\le \frac{1}{4}(|s|+|t|^{\frac{3}{2}})$. Similarly, when $|s|<\frac{1}{4\sqrt{M}}$ and $|t|<\frac{1}{64M}$, we can prove that
        $$
        |F_t(x,s,t)|\le 4M(4|t|^3+s^4) \le \frac{1}{4}(s^2+t^2).
        $$
        Hence, when $|(s,t)|\le \frac{1}{64M}$,
        \begin{eqnarray*}
        F_s(x,s,t)| &\le & \frac{1}{4}(|s|+|t|^{\frac{3}{2}}),\\
        F_t(x,s,t)| &\le & \frac{1}{4}(s^2+t^2).
        \end{eqnarray*}
        It is that $(C_1)$ holds.
 \item  When $s>1,F(x,s,s)=4s^4\ln(1+2s^4)\ge 4(s^2+s^3)$. So, $F$ satisfies $(C_2)$.
 \item  Let $\nu=4$ and $A=\frac{1}{4}$. For all $x\in V$, we have
        $$
             4F(x,s,t)-F_s(x,s,t)s-F_t(x,s,t)t
          =  -4M\frac{(s^4+t^4)^2}{(1+s^4+t^4)}
        \le \frac{1}{4}(s^2+|t|^3).
        $$
        So, $(C_3)$ holds.
 \item  Let $\beta_3=2$ and $K_3(x)\equiv 1$. For all $s\in \R$, $F(x,s,s)=2Ms^4\ln(1+2s^4)\ge -s^2$. So, $(C_4)$ holds.
\end{itemize}
Hence,  by Theorem 1.3, when $0<\lambda<\frac{\min\{\frac{1}{128M}\min\{\mu_0^{\frac{1}{2}},\mu_0^{\frac{1}{3}}\},1\}^2}{32\max\{\|e_1\|_{L^2(V)},\|e_2\|_{L^{\frac{3}{2}}(V)}\}}$,
system (\ref{eaxmple 2}) has one nontrivial solution of positive energy and another nontrivial solution of negative energy.

\vskip2mm
\section{Conclusion}
\vskip0mm
\noindent
The existence of  nontrivial solutions for system (\ref{eq1}) is investigated when the nonlinear term $F$ satisfies the sub-$(p,q)$-linear condition or super-$(p,q)$-linear condition,
which generalize some results in \cite{Grigor'yan 2017} in some sense. We present the concrete ranges of the parameter $\lambda_1$ and $\lambda_2$. For the sub-$(p,q)$-linear case, we
furthermore obtain a necessary condition for the existence of the non-semi-trivial solutions, and for the super-$(p,q)$-linear case, we present a weaker assumption of $F$ than the well-known (AR)-condition. However,
 we do not investigate the existence of the non-semi-trivial solutions. A possible method to solve the problem can be referred to \cite{Chang1} and \cite{Chang2} and we shall try to do it in future works.

\vskip2mm
\section{Appendix A}
  \setcounter{equation}{0}
In this section, we present some conclusions about $W_h^{1,s}(V)$ and $\varphi_{\lambda}$.
\qed
\vskip2mm
 \noindent
{\bf Appendix A.1.} {\it $W_h^{1,s}(V)$ is uniformly convex for all $s>1$.}
\vskip0mm
\noindent
{\bf Proof.}  Since $L^s(V)$ is uniformly convex for all $s>1$, by using Theorem 8 in \cite{xingjiasheng 2001}, we have $E$ is uniformly convex, where $E=L^s(V) \times L^s(V)$ with $\|(u,v)\|_E=\left(\|u\|_{L^s(V)}^s+\|v\|_{L^s(V)}^s\right)^{\frac{1}{s}}$.
Define $T:W_h^{1,s}(V)\to E$ by
$$
T(u(x))=\left(\nabla u(x),h(x)^{\frac{1}{s}}u(x)\right),
$$
where $h(x)\ge h_0>0$. Then
\begin{eqnarray*}
\|T(u)\|_E & = & \left(\|\nabla u\|_{L^s(V)}^s+\|h^{\frac{1}{s}}u)\|_{L^s(V)}^s\right)^{\frac{1}{s}}\\
           & = & \left(\int_V |\nabla u(x)|^s+h(x)|u(x)|^s d\mu\right)^{\frac{1}{s}}\\
           & = & \|u\|_{W_h^{1,s}(V)}.
\end{eqnarray*}
So, $T$ is an isometry. Hence, $W_h^{1,s}(V)$ is uniformly convex.\qed

\vskip2mm
\noindent
{\bf Appendix A.2.} {\it If $F(x,s,t)$satisfies $(F_0)$, then $\varphi_{\lambda}\in C^1(W,\R)$, and
\begin{eqnarray*}
       \langle\varphi'(u,v),(\phi_1,\phi_2)\rangle
 & = & \int_V\left[|\nabla u|^{p-2}\Gamma(u,\phi_1)+h_1|u|^{p-2}u\phi_1-F_u(x,u,v)\phi_1-\lambda_1 e_1\phi_1\right]d\mu\\
 &   & +\int_V\left[|\nabla v|^{q-2}\Gamma(v,\phi_2)+h_2|v|^{q-2}v\phi_2-F_v(x,u,v)\phi_2-\lambda_2 e_2\phi_2\right]d\mu.
\end{eqnarray*}}
\vskip0mm
\noindent
{\bf Proof.} Let
\begin{eqnarray*}
       \notag G(x,u,v)
 & = & \frac{1}{p}\left(|\nabla u|^p+h_1|u|^p\right)\mu(x)+\frac{1}{q}\left(|\nabla v|^q+h_2|v|^q\right)\mu(x)\\
 &   & -F(x,u,v)\mu(x)-\lambda_1e_1u \mu(x)-\lambda_2e_2v\mu(x).
\end{eqnarray*}
Then $\sum_{x\in V}G(x,u,v)=\varphi_{\lambda}(u,v)$. For any given $(\phi_1,\phi_2)\in W$ and $\theta \in [-1,1]$, we have
\begin{eqnarray*}
 G_x(\theta)
 & \triangleq & G(x,u+\theta \phi_1,v+\theta \phi_2)\\
 &     =      & \frac{1}{p}\left(|\nabla(u+\theta\phi_1)|^p+h_1|u+\theta\phi_1|^p\right)\mu(x)\\
 &            & +\frac{1}{q}\left(|\nabla(v+\theta\phi_2)|^q+h_2|v+\theta \phi_2|^q\right)\mu(x)\\
 &            & -F(x,u+\theta \phi_1,v+\theta \phi_2)\mu(x)-\lambda_1e_1(u+\theta \phi_1) \mu(x)-\lambda_2e_2(v+\theta \phi_2)\mu(x).
\end{eqnarray*}
Hence, by $(F_0)$ and (\ref{G2}),
we have
\begin{eqnarray}\label{A31}
\notag &      & G_x(\theta)\\
\notag & \le  & \frac{2^{p-1}}{p}\left(|\nabla u|^p+h_1|u|^p+|\nabla \phi_1|^p+h_1|\phi_1|^p\right)\mu(x)+\frac{2^{q-1}}{q}\left(|\nabla v|^q+h_2|v|^q+|\nabla \phi_2|^q+h_2|\phi_2|^q\right)\mu(x)\\
\notag &      & +a\left(|(u+\theta \phi_1,v+\theta \phi_2\right)|)b(x)\mu(x)+\lambda_1 |e_1(u+\theta \phi_1)|\mu(x)+\lambda_2 |e_2(v+\theta \phi_2)|\mu(x)\\
\notag & \le  & \frac{2^{p-1}}{p}\left(|\nabla u|^p+h_1|u|^p+|\nabla \phi_1|^p+h_1|\phi_1|^p\right)\mu(x)+\frac{2^{q-1}}{q}\left(|\nabla v|^q+h_2|v|^q+|\nabla \phi_2|^q+h_2|\phi_2|^q\right)\mu(x)\\
\notag &      & +\max_{|(s,t)| \le \|u\|_{\infty}+\|v\|_{\infty}+\|\phi_1\|_{\infty}+\|\phi_2\|_{\infty}}a(|(s,t)|)b(x)\mu(x)+\lambda_1 |e_1(u+ \phi_1)|\mu(x)
       +\lambda_2 |e_2(v+ \phi_2)|\mu(x).\\
\end{eqnarray}
Since $u,\phi_1 \in W_h^{1,p}(V)$, $v,\phi_2 \in W_h^{1,q}(V)$, $a\in L^{\infty}(V)$, $b\in L^1(V)$, $e_1 \in L^{\frac{p}{p-1}}(V)$, $e_2 \in L^{\frac{q}{q-1}}(V)$, then $\sum_{x \in V}G_x(\theta)$ is convergence
for all $\theta\in [-1,1]$.
Moreover,
\begin{eqnarray*}
&   & \frac{1}{p}\frac{\partial}{\partial \theta} |\nabla(u+\theta\phi_1)|^p\\
& = & \frac{1}{p}\frac{\partial}{\partial \theta} (|\nabla(u+\theta\phi_1)|^2)^{\frac{p}{2}}\\
& = & \frac{1}{2}|\nabla(u+\theta\phi_1)|^{p-2}\frac{\partial}{\partial \theta} \Gamma(u+\theta\phi_1,u+\theta\phi_1)\\
& = & \frac{1}{2}|\nabla(u+\theta\phi_1)|^{p-2}\frac{\partial}{\partial \theta} \left(\frac{1}{2\mu(x)}\sum\limits_{y\thicksim x}w_{xy}[(u(y)+\theta\phi_1(y))-(u(x)+\theta\phi_1(x))]^2\right)\\
& = & |\nabla(u+\theta\phi_1)|^{p-2}\frac{1}{2\mu(x)}\sum\limits_{y\thicksim x}w_{xy}[(u(y)+\theta\phi_1(y))-(u(x)+\theta\phi_1(x))](\phi_1(y)-\phi_1(x))\\
& = & |\nabla(u+\theta\phi_1)|^{p-2}\frac{1}{2\mu(x)}\sum\limits_{y\thicksim x}w_{xy}\left[(u(y)-u(x))(\phi_1(y)-\phi_1(x))+\theta(\phi_1(y)-\phi_1(x))^2\right]\\
& = & |\nabla(u+\theta\phi_1)|^{p-2}\left(\Gamma(u,\phi_1)+\Gamma(\theta\phi_1,\phi_1)\right)\\
& = & |\nabla(u+\theta\phi_1)|^{p-2}\Gamma(u+\theta\phi_1,\phi_1).
\end{eqnarray*}
Then
\begin{eqnarray}\label{A32}
\notag &   & \frac{\partial G_x(\theta)}{\partial \theta}\\
\notag & = & |\nabla (u+\theta\phi_1)|^{p-2}\Gamma(u+\theta\phi_1,\phi_1)\mu(x)\\
\notag &   & +|\nabla( v+\theta \phi_2)|^{q-2}\Gamma(v+\theta\phi_2, \phi_2)\mu(x)\\
\notag &   & +h_1|u+\theta\phi_1|^{p-2}(u+\theta\phi_1)\phi_1\mu(x)+h_2|v+\theta\phi_2|^{q-2}(v+\theta\phi_2)\phi_2\mu(x)\\
\notag &   & -F_{u+\theta \phi_1}(x,u+\theta \phi_1,v+\theta \phi_2)\phi_1\mu(x)-F_{v+\theta \phi_2}(x,u+\theta \phi_1,v+\theta \phi_2)\phi_2\mu(x)\\
       &   & -\lambda_1e_1 \phi_1 \mu(x)-\lambda_2e_2 \phi_2\mu(x).
\end{eqnarray}
Since $F(x,s,t)$ is continuously differentiable in $(s,t)\in \R^2$ for all $x\in V$, it is easy to obtain that $\frac{\partial G_x(\theta)}{\partial \theta}$ is continuous in $[-1,1]$.
By $(F_0)$, we have
\begin{eqnarray}\label{A33}
\notag &     & \frac{\partial G_x(\theta)}{\partial \theta}\\
\notag & \le & |\nabla (u+\theta\phi_1)|^{p-1}|\nabla \phi_1|\mu(x)+\left(h_1^{\frac{1}{p}}|u+\theta \phi_1|\right)^{p-1}h_1^{\frac{1}{p}}|\phi_1|\mu(x)\\
\notag &     & +|\nabla (v+\theta\phi_2)|^{q-1}|\nabla \phi_2|\mu(x)+\left(h_2^{\frac{1}{q}}|v+\theta \phi_2|\right)^{q-1}h_2^{\frac{1}{q}}|\phi_2|\mu(x)\\
\notag &     & +\left(|F_u|\phi_1+|F_v|\phi_2+\theta\lambda_1e_1 \phi_1 +\theta\lambda_2e_2\phi_2\right)\mu(x)\\
\notag & \le & \left(|\nabla (u+\theta\phi_1)|^{p}+h_1|u+\theta \phi_1|^p\right)^{\frac{p-1}{p}}\left(|\nabla \phi_1|^p+h_1|\phi_1|^p\right)^{\frac{1}{p}}\mu(x)\\
\notag &     & +\left(|\nabla (v+\theta \phi_2)|^{q}+h_2|v+\theta \phi_2|^q\right)^{\frac{q-1}{q}}\left(|\nabla \phi_2|^q+h_2|\phi_2|^q\right)^{\frac{1}{q}}\mu(x)\\
\notag &     & +\max_{|(s,t)| \le \|u\|_{\infty}+\|v\|_{\infty}}a|(s,t)|b(x)(\phi_1+\phi_2)\mu(x)+\lambda_1|e_1\phi_1|\mu(x)+\lambda_2|e_2\phi_2|\mu(x)\\
\notag & \le & 2^{\frac{(p-1)^2}{p}}\left(|\nabla u|^p+|\nabla \phi_1|^{p}+h_1|u|^p+|\phi_1|^p\right)^{\frac{p-1}{p}}\left(|\nabla \phi_1|^p+h_1|\phi_1|^p\right)^{\frac{1}{p}}\mu(x)\\
\notag &     & +2^{\frac{(q-1)^2}{q}}\left(|\nabla v|^q+|\nabla \phi_2|^{q}+h_2|v|^q+|\phi_2|^q\right)^{\frac{q-1}{q}}\left(|\nabla \phi_2|^q+h_2|\phi_2|^q\right)^{\frac{1}{q}}\mu(x)\\
\notag &     & +\max_{|(s,t)| \le \|u\|_{\infty}+\|v\|_{\infty}}a|(s,t)|b(x)(\phi_1+\phi_2)\mu(x)+\lambda_1|e_1\phi_1|\mu(x)+\lambda_2|e_2\phi_2|\mu(x)\\
\notag & \le & 2^{\frac{(p-1)^2}{p}}\left((|\nabla u|^p+h_1|u|^p)^{\frac{p-1}{p}}+(|\nabla \phi_1|^p+h_1|\phi_1|^p)^{\frac{p-1}{p}}\right)
               \left(|\nabla \phi_1|^p+h_1|\phi_1|^p\right)^{\frac{1}{p}}\mu(x)\\
\notag &     & +2^{\frac{(q-1)^2}{q}}\left((|\nabla v|^q+h_2|v|^q)^{\frac{q-1}{q}}+(|\nabla \phi_2|^q+h_2|\phi_2|^q)^{\frac{q-1}{q}}\right)
               \left(|\nabla \phi_2|^q+h_2|\phi_2|^q\right)^{\frac{1}{q}}\mu(x)\\
       &     & +\max_{|(s,t)| \le \|u\|_{\infty}+\|v\|_{\infty}}a|(s,t)|b(x)(\phi_1+\phi_2)\mu(x)+\lambda_1|e_1\phi_1|\mu(x)+\lambda_2|e_2\phi_2|\mu(x).
\end{eqnarray}
Moreover, we have
\begin{eqnarray}\label{A35}
\notag &     &  \sum_{x \in V}\frac{\partial G_x(\theta)}{\partial \theta}\\
\notag & \le & 2^{\frac{(p-1)^2}{p}}\left(\sum_{x\in V}\left((|\nabla u|^p+h_1|u|^p)^{\frac{p-1}{p}}+(|\nabla \phi_1|^p+h_1|\phi_1|^p)^{\frac{p-1}{p}}\right)^
               {\frac{p}{p-1}}\right)^{\frac{p-1}{p}}\left(\sum_{x\in V}|\nabla \phi_1|^p+h_1|\phi_1|^p\right)^{\frac{1}{p}}\mu(x)\\
\notag &     & +2^{\frac{(q-1)^2}{q}}\left(\sum_{x\in V}\left((|\nabla v|^q+h_2|v|^q)^{\frac{q-1}{q}}+(|\nabla \phi_2|^q+h_2|\phi_2|^q)^{\frac{q-1}{q}}\right)^
               {\frac{q}{q-1}}\right)^{\frac{q-1}{q}}\left(\sum_{x\in V}|\nabla \phi_2|^q+h_2|\phi_2|^q\right)^{\frac{1}{q}}\mu(x)\\
\notag &     & +\max_{|(s,t)| \le \|u\|_{\infty}+\|v\|_{\infty}}a|(s,t)|\|b(x)\|_{L^1(V)}\|\phi_1+\phi_2\|_{\infty}+\lambda_1\|e_1\|_{L^{\frac{p}{p-1}}(V)}\|\phi_1\|_{L^p(V)}
               +\lambda_2\|e_2\|_{L^{\frac{q}{q-1}}(V)}\|\phi_2\|_{L^q(V)}\\
\notag & \le & 2^{\frac{(p-1)^2+1}{p}}\left(\sum_{x\in V}\left(|\nabla u|^p+h_1|u|^p+|\nabla \phi_1|^p+h_1|\phi_1|^p\right)\mu(x)\right)^{\frac{p-1}{p}}
               \left(\sum_{x\in V}(|\nabla \phi_1|^p+h_1|\phi_1|^p)\mu(x)\right)^{\frac{1}{p}}\\
\notag &     & +2^{\frac{(q-1)^2+1}{q}}\left(\sum_{x\in V}\left(|\nabla v|^q+h_2|v|^q+|\nabla \phi_2|^q+h_2|\phi_2|^q\right)\mu(x)\right)^{\frac{q-1}{q}}
               \left(\sum_{x\in V}(|\nabla \phi_2|^q+h_2|\phi_2|^q)\mu(x)\right)^{\frac{1}{q}}\\
\notag &     & +\max_{|(s,t)| \le \|u\|_{\infty}+\|v\|_{\infty}}a|(s,t)|\|b(x)\|_{L^1(V)}\|\phi_1+\phi_2\|_{\infty}+\lambda_1\|e_1\|_{L^{\frac{p}{p-1}}(V)}\|\phi_1\|_{L^p(V)}
               +\lambda_2\|e_2\|_{L^{\frac{q}{q-1}}(V)}\|\phi_2\|_{L^q(V)}\\
\notag &  =  & 2^{\frac{(p-1)^2+1}{p}}\left(\|u\|_{W^{1,p}(V)}^p+\|\phi_1\|_{W^{1,p}(V)}^p\right)^{\frac{p-1}{p}}\|\phi_1\|_{W^{1,p}(V)}
               +\lambda_1\|e_1\|_{L^{\frac{p}{p-1}}(V)}\|\phi_1\|_{L^p(V)}\\
\notag &     & 2^{\frac{(q-1)^2+1}{q}}\left(\|v\|_{W^{1,q}(V)}^q+\|\phi_2\|_{W^{1,q}(V)}^q\right)^{\frac{q-1}{q}}\|\phi_2\|_{W^{1,q}(V)}
               +\lambda_2\|e_2\|_{L^{\frac{q}{q-1}}(V)}\|\phi_2\|_{L^q(V)}\\
\notag &     & +\max_{|(s,t)| \le \|u\|_{\infty}+\|v\|_{\infty}}a|(s,t)|\|b(x)\|_{L^1(V)}\|\phi_1+\phi_2\|_{\infty}\\
\notag & \le & 2^{\frac{(p-1)^2+1}{p}}\left(\|u\|_{W^{1,p}(V)}^{p-1}+\|\phi_1\|_{W^{1,p}(V)}^{p-1}\right)\|\phi_1\|_{W^{1,p}(V)}
               +\lambda_1\|e_1\|_{L^{\frac{p}{p-1}}(V)}\|\phi_1\|_{L^p(V)}\\
\notag &     & + 2^{\frac{(q-1)^2+1}{q}}\left(\|v\|_{W^{1,q}(V)}^{q-1}+\|\phi_2\|_{W^{1,q}(V)}^{q-1}\right)\|\phi_2\|_{W^{1,q}(V)}
               +\lambda_2\|e_2\|_{L^{\frac{q}{q-1}}(V)}\|\phi_2\|_{L^q(V)}\\
       &     & +\max_{|(s,t)| \le \|u\|_{\infty}+\|v\|_{\infty}}a|(s,t)|\|b(x)\|_{L^1(V)}\|\phi_1+\phi_2\|_{\infty}.
\end{eqnarray}
So, we obtain that $\sum_{x \in V}\frac{\partial G_x(\theta)}{\partial \theta}$ is uniform convergence.
Let $H(\theta)=\sum_{x \in V}\frac{\partial G_x(\theta)}{\partial \theta}=\varphi_{\lambda}(u+\theta \phi_1,v+\theta \phi_2)$. Then by (\ref{A31}), (\ref{A32}), (\ref{A33}) and (\ref{A35}), we have
\begin{eqnarray}\label{A34}
\notag  H'(0)
        &  =  & \left(\sum_{x \in V}\frac{\partial G_x(\theta)}{\partial \theta}\right)\big|_{\theta=0}\\
\notag  &  =  &  \sum_{x \in V}\left[|\nabla u|^{p-2}\Gamma(u,\phi_1)+h_1|u|^{p-2}u\phi_1-F_u(x,u,v)\phi_1-\lambda_1 e_1\phi_1\right]\mu(x)\\
\notag  &     &  +\sum_{x \in V}\left[|\nabla v|^{q-2}\Gamma(v,\phi_2)+h_2|v|^{q-2}v\phi_2-F_v(x,u,v)\phi_2-\lambda_2 e_2\phi_2\right]\mu(x)\\
        &  =  &  \langle\varphi_{\lambda}'(u,v),(\phi_1,\phi_2)\rangle.
\end{eqnarray}
So, for any given $(\phi_1,\phi_2) \in W$, by (\ref{A34}) and (\ref{G1}), we have
\begin{eqnarray*}
&     & \langle\varphi_{\lambda}'(u,v),(\phi_1,\phi_2)\rangle\\
& \le & \int_V\left[|\nabla u|^{p-1}|\nabla\phi_1|+h_1|u|^{p-2}u\phi_1-F_u(x,u,v)\phi_1-\lambda_1 e_1\phi_1\right]\mu(x)\\
&     &  +\int_V\left[|\nabla v|^{q-1}|\nabla\phi_2|+h_2|v|^{q-2}v\phi_2-F_v(x,u,v)\phi_2-\lambda_2 e_2\phi_2\right]\mu(x)\\
& \le & \left(\int_V|\nabla u|^pd\mu\right)^{\frac{p-1}{p}}\left(\int_V|\nabla \phi_1|^pd\mu\right)^{\frac{1}{p}}+\left(\int_Vh_1|u|^pd\mu\right)^{\frac{p-1}{p}}\left(\int_Vh_1|\phi_1|^pd\mu\right)^{\frac{1}{p}}\\
&     & +\left(\int_V|\nabla v|^qd\mu\right)^{\frac{q-1}{q}}\left(\int_V|\nabla \phi_2|^qd\mu\right)^{\frac{1}{q}}+\left(\int_Vh_2|v|^qd\mu\right)^{\frac{q-1}{q}}\left(\int_Vh_2|\phi_2|^qd\mu\right)^{\frac{1}{q}}\\
&     & +\int_V\left(|F_u|\phi_1+|F_v|\phi_2+\lambda_1e_1\theta \phi_1 +\lambda_2e_2\theta \phi_2\right)d\mu\\
& \le & \|u\|_{W^{1,p}(V)}^{p-1}\|\phi_1\|_{W^{1,p}(V)}+\|v\|_{W^{1,q}(V)}^{q-1}\|\phi_2\|_{W^{1,q}(V)}\\
&     & +\int_Va(|(u,v)|)b(x)(\phi_1+\phi_2)d\mu+\int_V(\lambda_1e_1\phi_1 +\lambda_2e_2 \phi_2)d\mu\\
& \le & \|u\|_{W^{1,p}(V)}^{p-1}\|\phi_1\|_{W^{1,p}(V)}+\|v\|_{W^{1,q}(V)}^{q-1}\|\phi_2\|_{W^{1,q}(V)}\\
&     & +\max_{|(s,t)| \le \|u\|_{\infty}+\|v\|_{\infty}}a|(s,t)|\|b\|_{L^1(V)}(\|\phi_1\|_{\infty}+\|\phi_2\|_{\infty})\\
&     & +\lambda_1h_0^{-\frac{1}{p}}\|e_1\|_{L^{\frac{p}{p-1}}(V)}\|\phi_1\|_{W^{1,p}(V)}+\lambda_2h_0^{-\frac{1}{q}}\|e_2\|_{L^{\frac{q}{q-1}}(V)}\|\phi_2\|_{W^{1,q}(V)}\\
& \le & \max\{\|u\|_{W^{1,p}(V)}^{p-1},\|v\|_{W^{1,q}(V)}^{q-1}\}\|(\phi_1,\phi_2)\|_W\\
&     & +\max\{\frac{1}{(h_0\mu_0)^{\frac{1}{p}}},\frac{1}{(h_0\mu_0)^{\frac{1}{q}}}\}\max_{|(s,t)| \le \|u\|_{\infty}+\|v\|_{\infty}}a|(s,t)|\|b\|_{L^1(V)}\|(\phi_1,\phi_2)\|_W\\
&     & +\max\{\lambda_1h_0^{-\frac{1}{p}}\|e_1\|_{L^{\frac{p}{p-1}}(V)},\lambda_2h_0^{-\frac{1}{q}}\|e_2\|_{L^{\frac{q}{q-1}}(V)}\}\|(\phi_1,\phi_2)\|_W.
\end{eqnarray*}
Thus, $\varphi_{\lambda}'(u,v):W\to \R$ is bounded and linear operator, that is, $\varphi_{\lambda}'(u,v)\in W^*$ which is the dual space of $W$.
Define the mapping $\varphi_{\lambda}' :W\to W^*$ by
$$
\varphi_{\lambda}': (u,v) \to \varphi_{\lambda}'(u,v).
$$
Next, we prove that $\varphi_{\lambda}'$ is continuous in $W$. For any sequence $\{(u_k,v_k)\} \subset W$ with $(u_k,v_k)\rightarrow (u,v)$  in $W$ as $k\rightarrow \infty$, we have
\begin{eqnarray}\label{k2}
\int_V|\nabla (u_k-u)|^pd\mu \to 0,\ \int_V|(u_k-u)|^pd\mu \to 0,\  \int_V|\nabla (v_k-v)|^qd\mu \to 0,\ \int_V|(v_k-v)|^qd\mu \to 0
\end{eqnarray}
 and by (\ref{bb1}), we have
\begin{eqnarray}\label{k1}
u_k(x)\to u(x),\ \ v_k(x)\to v(x),\ \mbox{for all }x\in V, \ \mbox{as }k\to \infty.
\end{eqnarray}
Note that
\begin{eqnarray*}
&     & \langle\varphi_{\lambda}'(u,v)-\varphi_{\lambda}'(u_k,v_k),(\phi_1,\phi_2)\rangle\\
&  =  & \int_V\left[|\nabla u|^{p-2}\Gamma(u,\phi_1)-|\nabla u_k|^{p-2}\Gamma(u_k,\phi_1)+h_1(|u|^{p-2}u-|u_k|^{p-2}u_k)\phi_1\right]d\mu\\
&     & +\int_V\left[|\nabla v|^{q-2}\Gamma(v,\phi_2)-|\nabla v_k|^{q-2}\Gamma(v_k,\phi_2)+h_2(|v|^{q-2}v-|v_k|^{q-2}v_k)\phi_2\right]d\mu\\
&     & -\int_V(F_u(x,u,v)-F_{u_k}(x,u_k,v_k))\phi_1d\mu-\int_V(F_v(x,u,v)-F_{v_k}(x,u_k,v_k))\phi_2d\mu\\
&  =  & \int_V|\nabla u|^{p-2}\Gamma(u-u_k,\phi_1)d\mu+\int_V(|\nabla u|^{p-2}-|\nabla u_k|^{p-2})\Gamma(u_k,\phi_1)d\mu\\
&     & +\int_V h_1(|u|^{p-2}u-|u_k|^{p-2}u_k)\phi_1d\mu-\int_V(F_u(x,u,v)-F_{u_k}(x,u_k,v_k))\phi_1d\mu\\
&     & +\int_V|\nabla v|^{q-2}\Gamma(v-v_k,\phi_2)d\mu+\int_V(|\nabla v|^{q-2}-|\nabla v_k|^{q-2})\Gamma(v_k,\phi_2)d\mu\\
&     & +\int_V h_2(|v|^{q-2}v-|v_k|^{q-2}v_k)\phi_2d\mu-\int_V(F_v(x,u,v)-F_{v_k}(x,u_k,v_k))\phi_2d\mu\\
&  :=  & I+II.
\end{eqnarray*}
Firstly, we prove that
\begin{eqnarray*}
I
&  =  & \int_V|\nabla u|^{p-2}\Gamma(u-u_k,\phi_1)d\mu+\int_V(|\nabla u|^{p-2}-|\nabla u_k|^{p-2})\Gamma(u_k,\phi_1)d\mu\\
&     & +\int_V h_1(|u|^{p-2}u-|u_k|^{p-2}u_k)\phi_1d\mu-\int_V(F_u(x,u,v)-F_{u_k}(x,u_k,v_k))\phi_1d\mu\\
&\rightarrow& 0 \mbox{ as }k\rightarrow\infty.
\end{eqnarray*}
By using Lemma 5.12 in \cite{Meng}, we have
\begin{eqnarray}\label{C11}
&     & \int_V h_1(|u|^{p-2}u-|u_k|^{p-2}u_k)\phi_1d\mu\nonumber\\
& \le &  \left(\int_V h_1(|u|^{p-2}u-|u_k|^{p-2}u_k)^{\frac{p}{p-1}}d\mu\right)^{\frac{p-1}{p}}\left(\int_V h_1|\phi_1|^pd\mu\right)^{\frac{1}{p}}\nonumber\\
&\rightarrow& 0 \mbox{ as }k\rightarrow\infty.
\end{eqnarray}
Similarly, by $(F_0)$, Lebesgue dominated convergence theorem, (\ref{k1}) and the continuity of $F_u$, we also have
\begin{eqnarray}\label{C12}
&     & \int_V(F_u(x,u,v)-F_{u_k}(x,u_k,v_k))\phi_1d\mu\nonumber\\
& \le & \left(\int_V(F_u(x,u,v)-F_{u_k}(x,u_k,v_k))^{\frac{p}{p-1}}d\mu\right)^{\frac{p-1}{p}}\left(\int_V |\phi_1|^pd\mu\right)^{\frac{1}{p}}\nonumber\\
&\rightarrow& 0 \mbox{ as }k\rightarrow\infty.
\end{eqnarray}
Moreover, by H\"older inequality and (\ref{k2}), we get
\begin{eqnarray}\label{C13}
&     & \int_V|\nabla u|^{p-2}\Gamma(u-u_k,\phi_1)d\mu\nonumber\\
& \le & \int_V|\nabla u|^{p-2}|\nabla(u-u_k)|\cdot|\nabla\phi_1|d\mu\nonumber\\
& \le & \left(\int_V|\nabla(u-u_k)|^pd\mu\right)^{\frac{1}{p}}\left(\int_V|\nabla u|^{\frac{(p-2)p}{p-1}}\cdot|\nabla\phi_1|^{\frac{p}{p-1}}d\mu\right)^{\frac{p-1}{p}}\nonumber\\
& \le & \left(\int_V|\nabla(u-u_k)|^pd\mu\right)^{\frac{1}{p}}\left(\int_V|\nabla u|^pd\mu\right)^{\frac{p-2}{p}}\left(\int_V|\nabla\phi_1|^pd\mu\right)^{\frac{1}{p}}\nonumber\\
&\rightarrow& 0 \mbox{ as }k\rightarrow\infty,
\end{eqnarray}
and
\begin{eqnarray}\label{C14}
&     & \int_V(|\nabla u|^{p-2}-|\nabla u_k|^{p-2})\Gamma(u_k,\phi_1)d\mu\nonumber\\
& \le & \int_V\big||\nabla u|^{p-2}-|\nabla u_k|^{p-2}\big||\nabla u_k|\cdot|\nabla\phi_1|d\mu\nonumber\\
&  =  & \int_V\big||\nabla u|^{p-2}|\nabla u_k|-|\nabla u_k|^{p-1}\big||\nabla\phi_1|d\mu\nonumber\\
&  =  & \int_V\big||\nabla u|^{p-2}\left[|\nabla u|+(|\nabla u_k|-|\nabla u|)\right]-|\nabla u_k|^{p-1}\big||\nabla\phi_1|d\mu\nonumber\\
&  =  & \int_V\big||\nabla u|^{p-1}-|\nabla u_k|^{p-1}\big||\nabla\phi_1|d\mu+\int_V|\nabla u|^{p-2}\big||\nabla u_k|-|\nabla u|\big||\nabla\phi_1|d\mu\nonumber\\
& \le & \left(\int_V\big||\nabla u|^{p-1}-|\nabla u_k|^{p-1}\big|^{\frac{p}{p-1}}d\mu\right)^{\frac{p-1}{p}}\left(\int_V|\nabla\phi_1|^pd\mu\right)^{\frac{1}{p}}\nonumber\\
&     & +\left(\int_V|\nabla u|^{\frac{(p-2)p}{p-1}}\big||\nabla u_k|-|\nabla u|\big|^{\frac{p}{p-1}}d\mu\right)^{\frac{p-1}{p}}\left(\int_V|\nabla\phi_1|^pd\mu\right)^{\frac{1}{p}}\nonumber\\
& \le & (p-1)\left(\int_V|\nabla\phi_1|^pd\mu\right)^{\frac{1}{p}}\left(\int_V\big||\nabla u|-|\nabla u_k|\big|^{\frac{p}{p-1}}\left(|\nabla u_k|^{p-2}+|\nabla u|^{p-2}\right)^{\frac{p}{p-1}}d\mu\right)^{\frac{p-1}{p}}\nonumber\\
&     & +\left(\int_V|\nabla\phi_1|^pd\mu\right)^{\frac{1}{p}}\left(\int_V|\nabla u|^pd\mu\right)^{\frac{p-2}{p}}\left(\int_V\big||\nabla u_k|-|\nabla u|\big|^pd\mu\right)^{\frac{1}{p}}\nonumber\\
& \le & (p-1)\left(\int_V|\nabla\phi_1|^pd\mu\right)^{\frac{1}{p}}\left(\int_V|\nabla (u-u_k)|^pd\mu\right)^{\frac{1}{p}}\left(\int_V(|\nabla u_k|^{p-2}+|\nabla u|^{p-2})^{\frac{p}{p-2}}d\mu\right)^{\frac{p-2}{p}}\nonumber\\
&     & +\left(\int_V|\nabla\phi_1|^pd\mu\right)^{\frac{1}{p}}\left(\int_V|\nabla u|^pd\mu\right)^{\frac{p-2}{p}}\left(\int_V|\nabla (u_k-u)|^pd\mu\right)^{\frac{1}{p}}\nonumber\\
& \le & 2^{\frac{2}{p}}(p-1)\left(\int_V|\nabla\phi_1|^pd\mu\right)^{\frac{1}{p}}\left(\int_V|\nabla (u-u_k)|^pd\mu\right)^{\frac{1}{p}}\left(\int_V(|\nabla u_k|^{p}+|\nabla u|^{p})d\mu\right)^{\frac{p-2}{p}}\nonumber\\
&     & +\left(\int_V|\nabla\phi_1|^pd\mu\right)^{\frac{1}{p}}\left(\int_V|\nabla u|^pd\mu\right)^{\frac{p-2}{p}}\left(\int_V|\nabla (u_k-u)|^pd\mu\right)^{\frac{1}{p}}\nonumber\\
&\rightarrow& 0 \mbox{ as }k\rightarrow\infty.
\end{eqnarray}
So, by (\ref{C11})-(\ref{C14}), we obtain that
$$
I\rightarrow 0 \mbox{ as }k\rightarrow\infty.
$$
Similarly, we can prove that
$$
II\rightarrow 0 \mbox{ as }k\rightarrow\infty.
$$
Hence,
$$
\langle\varphi_{\lambda}'(u,v)-\varphi_{\lambda}'(u_k,v_k),(\phi_1,\phi_2)\rangle\rightarrow 0 \mbox{ as }k\rightarrow\infty.
$$
Then $\varphi_{\lambda}'$ is continuous.\qed

\vskip3mm
 \noindent
\noindent{\bf Funding information}

\noindent
This project is supported by Yunnan Ten Thousand Talents Plan Young \& Elite Talents Project and supported by Yunnan Fundamental Research Projects (grant No: 202301AT070465).

\vskip3mm
 \noindent
\noindent{\bf Author contribution}
\par
\noindent
 These authors contributed equally to this work.

\vskip3mm
 \noindent
\noindent{\bf Conflict of interest}

\noindent
The authors state no conflict of interest.

\vskip3mm
 \noindent
\noindent{\bf Data availability statement }

\noindent
Not available.
\vskip2mm
\renewcommand\refname{References}
{}
\end{document}